\documentclass[10pt,a4paper]{article}
\usepackage{bbm}

\usepackage[leqno]{amsmath}
\usepackage{amsfonts}
\usepackage{graphicx}
\usepackage{amsmath}
\usepackage{amssymb}
\usepackage{latexsym}
\usepackage{amsmath, amsfonts,amssymb, amsthm, euscript,makeidx,color,mathrsfs}
\usepackage{enumerate}

\usepackage[colorlinks,linkcolor=blue,anchorcolor=green,citecolor=red]{hyperref}

\def\sqr#1#2{{\vcenter{\vbox{\hrule height.#2pt
              \hbox{\vrule width.#2pt height#1pt \kern#1pt \vrule width.#2pt}
              \hrule height.#2pt}}}}
\def\signed #1{{\unskip\nobreak\hfil\penalty50
              \hskip2em\hbox{}\nobreak\hfil#1
              \parfillskip=0pt \finalhyphendemerits=0 \par}}
\def\endpf{\signed {$\sqr69$}}

\def\5n{\negthinspace \negthinspace \negthinspace \negthinspace \negthinspace }
\def\4n{\negthinspace \negthinspace \negthinspace \negthinspace }
\def\3n{\negthinspace \negthinspace \negthinspace }
\def\2n{\negthinspace \negthinspace }
\def\1n{\negthinspace }

\def\dbR{\mathbb{R}}

\def\sE{\mathscr{E}}

\def\sT{\mathscr{T}}
\def\sU{\mathscr{U}}

\def\={\buildrel \triangle \over =}

\def\ds{\displaystyle}

\def\ns{\noalign{\ss}}
\def\a{\alpha}
\def\b{\beta}
\def\g{\gamma}
\def\d{\delta}
\def\e{\varepsilon}

\def\k{\kappa}

\def\si{\sigma}
\def\t{\tau}
\def\f{\varphi}
\def\th{\theta}
\def\o{\omega}

\def\G{\Gamma}
\def\D{\Delta}

\def\F{\Phi}

\def\cA{{\cal A}}

\def\cU{{\cal U}}

\def\ss{\smallskip}
\def\ms{\medskip}
\def\bs{\bigskip}
\def\q{\quad}
\def\qq{\qquad}
\def\hb{\hbox}

\def\ua{\mathop{\uparrow}}

\def\Ra{\mathop{\Rightarrow}}

\def\esssup{\mathop{\rm esssup}}

\def\h{\widehat}
\def\wt{\widetilde}

\def\cd{\cdot}
\def\cds{\cdots}

\def\ae{\hbox{\rm a.e.}}

\def\deq{\triangleq}
\def\les{\leqslant}
\def\ges{\geqslant}

\def\({\Big (}
\def\){\Big )}
\def\[{\Big[}
\def\]{\Big]}

\def\bde{\begin{definition}\label}
\def\ede{\end{definition}}
\def\be{\begin{equation}}
\def\bel{\begin{equation}\label}
\def\ee{\end{equation}}
\def\bt{\begin{theorem}\label}
\def\et{\end{theorem}}
\def\bc{\begin{corollary}\label}
\def\ec{\end{corollary}}
\def\bl{\begin{lemma}\label}
\def\el{\end{lemma}}
\def\bp{\begin{proposition}\label}
\def\ep{\end{proposition}}
\def\bas{\begin{assumption}\label}
\def\eas{\end{assumption}}
\def\br{\begin{remark}\label}
\def\er{\end{remark}}
\def\bex{\begin{example}\label}
\def\ex{\end{example}}
\def\ba{\begin{array}}
\def\ea{\end{array}}
\def\ed{\end{document}}

\def\square#1{\vbox{\hrule\hbox{\vrule height#1%
     \kern#1\vrule}\hrule}}
\def\rectangle#1#2{\vbox{\hrule\hbox{\vrule height#1%
     \kern#2\vrule}\hrule}}

\font\tenbb=msbm10 \font\sevenbb=msbm7 \font\fivebb=msbm5

\newfam\bbfam
\scriptscriptfont\bbfam=\fivebb \textfont\bbfam=\tenbb
\scriptfont\bbfam=\sevenbb

\newtheorem{theorem}{\hskip 1.3em Theorem}[section]
\newtheorem{definition}[theorem]{\hskip 1.3em Definition}
\newtheorem{proposition}[theorem]{\hskip 1.3em Proposition}
\newtheorem{corollary}[theorem]{\hskip 1.3em Corollary}
\newtheorem{lemma}[theorem]{\hskip 1.3em Lemma}
\newtheorem{remark}[theorem]{\hskip 1.3em Remark}
\newtheorem{example}[theorem]{\hskip 1.3em Example}

\newtheorem{assumption}[theorem]{\hskip 1.3em Assumption}

\makeatletter
   
   \@addtoreset{equation}{section}
\makeatother

\begin{document}

\title{\bf Controlled Singular Volterra Integral Equations \\ and Pontryagin Maximum Principle \footnote{This work was partially supported by the National Natural Science Foundation of China under grant 11471070, and NSF Grant DMS-1406776.}}
\date{}
\author{Ping Lin\footnote{  School of
Mathematics and Statistics, Northeast Normal University, Changchun
130024, China. E-mail address:
linp258@nenu.edu.cn.}~~~and~~Jiongmin Yong \footnote{Department of Mathematics, University of Central Florida, Orlando, FL 32816, USA.  E-mail address: jiongmin.yong@ucf.edu.}}

\ms

\maketitle

{\bf Abstract.}
This paper is concerned with a class of controlled singular Volterra integral equations, which could be used to describe problems involving memories. The well-known fractional order ordinary differential equations of the Riemann--Liouville or Caputo types are strictly special cases of the   equations studied in this paper. Well-posedness and some regularity results in proper spaces are established for such kind of questions. For the associated optimal control problem, by using a Liapounoff's type theorem and the spike variation technique, we establish a Pontryagin's type maximum principle for optimal controls. Different from the existing literature, our method enables us to deal with the problem without assuming regularity conditions on the controls, the convexity condition on the control domain, and some additional unnecessary conditions on the nonlinear terms of the integral equation and the cost functional.

\bs


\ms


\section{Introduction.}

In this paper, we consider the following controlled Volterra integral equation:
\bel{1.1}y(t)=\eta(t)+\int_0^tf(t,s,y(s),u(s))ds,\qq t\in[0,T].\ee
We call the above the {\it state equation}, where $\eta(\cd)$ and $f(\cd\,,\cd\,,\cd\,,\cd)$ are given maps, called the {\it free term} and the {\it generator} of the state equation, respectively, $y(\cd)$ is called the {\it state trajectory} taking values in the Euclidean space $\dbR^n$, and $u(\cd)$ is called the {\it control} taking values in some separable metric space $U$. To measure the performance of the control, we introduce the cost functional
\bel{1.2}J(u(\cd))=\int_0^Tg(t,y(t),u(t))dt+\sum_{j=1}^m h^j(y(t_j)),\ee
with the two terms on the right hand representing the {\it running cost} and the {\it specific instant costs} (at $0\les t_1<t_2<\cds<t_m\les T$), respectively.

\ms

Equations like (\ref{1.1}) can be used to describe some dynamics involving memories. In the classical situations of optimal control for Volterra integral equations, people usually assume that the map $f(\cd\,,\cd\,,\cd\,,\cd)$ is continuous, together with some further smoothness/differentiability conditions. Relevant works can be traced back to those by Vinokurov in the later 1960s \cite{Vinokurov 1969}, followed by the works of Angell \cite{Angell 1976}, Kamien-Muller \cite{Kamien-Muller 1976}, Medhin \cite{Medhin 1986}, Carlson \cite{Carlson 1987}, Burnap-Kazemi \cite{Burnap-Kazemi 1999}, and some recent works by de la Vega \cite{de la Vega 2006}, Belbas
\cite{Belbas 2007,Belbas 2008}, and Bonnans--de la Vega--Dupuis \cite{Bonnans-de la Vega-Dupuis 2013}. On the other hand, in the past several decades, fractional differential equations have attracted quite a few researchers' attention due to some very interesting applications in physics, chemistry, engineering, population dynamics, finance and other sciences; See Oldham--Spanier \cite{Oldham-Spanier 1974} for some early examples of diffusion processes, Torvik--Bagley \cite{Torvik-Bagley 1984}, Caputo \cite{Caputo 1967}, and Caputo--Mainardi \cite{Caputo--Mainardi 1971} for modeling of the mechanical properties of materials, Benson \cite{Benson 1998} for the advection and the dispersion of solutes in natural porous or fractured media, Chern \cite{Chern 1993}, Diethelm--Freed \cite{Diethelm-Freed 1999} for the modeling  behavior of viscoelastic and viscoplastic materials under external influences,
Scalas--Gorenflo--Mainardi \cite{Scalas-Gorenflo-Mainardi 2004} for the mathematical models in finance, Das--Gupta \cite{Das-Gupta 2011}, Demirci--Unal--\"Ozalp \cite{Demirci-Unal-Ozalp 2011}, Arafa--Rida--Khalil \cite{Arafa-Rida-Khalil 2012}, Diethelm \cite{Diethelm 2013} for some population and epidemic models, and Metzler et al. \cite{Metzler-Schick-Kilian-Nonnenmacher 1995} for the relaxation in filled polymer networks. An extensive survey on fractional differential equations can be found in the book by Kilbas--Srivastava--Trujillo \cite{Kilbas-Srivastava-Trujillo 2006}. In the recent years, optimal control problems have been studied for fractional differential equations by a number of authors. We mention the works of Agrawal \cite{Agrawal 2004, Agrawal 2008}, Agrawal--Defterli--Baleanu \cite{Agrawal-Defterli-Baleanu 2010},  Bourdin \cite{Bourdin 2012}, Frederico--Torres \cite{Frederico-Torres 2008}, Hasan--Tangpong--Agrawal \cite{Hasan-Tangpong-Agrawal} and Kamocki \cite{Kamocki 2014,Kamocki 2014b}.

\ms

The most popular fractional differential equations are those in the sense of Riemann--Liouville or in the sense of Caputo (See Section 3 for some details). It turns out that these equations (of the order no more than 1, for scalar functions) are equivalent to Volterra integral equations with the generator being singular along $s=t$, and the free term $\eta(\cd)$ being possibly discontinuous (blowing up) at $t=0$. More precisely, the corresponding controlled state equation of form (\ref{1.1}) could have the feature that
\bel{1.3*}\eta(t)={c\over t^{1-\a}}\ (\mbox{or}\ c),\qq f(t,s,y,u)={\wt f(s,y,u)\over(t-s)^{1-\a}},\qq0\les s<t\les T,\q\forall(y,u),\ee
for some map $\wt f(\cd\,,\cd\,,\cd)$ and constants $\a\in(0,1)$, $c\in\dbR$. Such kind of singularity makes the optimal control problems for fractional differential equations different from the classical optimal control problems for Volterra integral equations as in the above-mentioned literature.

\ms

The purpose of this paper is to study an optimal control problem with the state equation (\ref{1.1}) allowing $(t,s)\mapsto f(t,s,y,u)$ to have some singularity along $t=s$ and allowing the free term $\eta(\cd)$ to be (unboundedly) discontinuous. We point out that our state equation (\ref{1.1}) could cover a much wider class of dynamic systems with various type memories than the ones described by fractional differential equations (with the conditions like (\ref{1.3*})). Let us make a little more comments on our state equation. Since the free term $\eta(\cd)$ is allowed to have some singularities, a natural class for $\eta(\cd)$ should be $L^p$ functions. Then we expect, under suitable conditions, the state trajectory $y(\cd)$ will also be a function in the same class. On the other hand, in the cost functional, we need $y(t_j)$ to be defined. Therefore, we need to have certain continuity of the state trajectory. Then it is necessary to narrow the $L^p$ space by adding certain continuity. This will lead to some difficulties in establishing the well-posedness of the state equation in the correct class of functions that the state trajectories will belong to. To overcome the difficulty, we introduce certain weighted function spaces, and extend some classical results, such as Gronwall's inequality, etc. to the form that will make our procedure works.

\ms

The rest of the paper is organized as follows. In Section 2, necessary preliminaries will be presented. Some results are interesting by themselves. Well-posedness of the state equation, together with the continuity of the solutions, will be established in Section 3. Section 4 is devoted to a proof of Pontryagin's type maximum principle for our optimal control problem of singular integral equations. As a special case, the maximum principles for fractional differential equations in the sense of Riemann--Liouville, and Caputo, will be briefly described. Some concluding remarks will be collected in Section 5.

\section{Preliminary}

In this section, we will present some preliminary results which will be useful later. First of all, let $T>0$ be a fixed time horizon. We introduce the following spaces:
$$\ba{ll}
\ns\ds L^p(0,T;\dbR^n)=\Big\{\f:[0,T]\to\dbR^n\bigm|\|\f(\cd)\|_{L^p(0,T;\dbR^n)}\equiv\(\int_0^T
|\f(t)|^pdt\)^{1\over p}<\infty\Big\},\q1\les p<\infty,\\
\ns\ds L^\infty(0,T;\dbR^n)=\Big\{\f:[0,T]\to\dbR^n\bigm|\|\f(\cd)\|_{L^\infty(0,T;\dbR^n)}\equiv
\esssup_{t\in[0,T]}|\f(t)|<\infty\Big\},\\
\ns\ds C([0,T];\dbR^n)=\Big\{\f:[0,T]\to\dbR^n\bigm|\f(\cd)\hb{ is continuous}\Big\}.\ea$$
We denote
$$\ba{ll}
\ns\ds\|\f(\cd)\|_p=\|\f(\cd)\|_{L^p(0,T;\dbR^n)},\qq\forall\f(\cd)\in L^p(0,T;\dbR^n),\qq p\in[1,\infty],\\
\ns\ds\|\f(\cd)\|_C=\max_{t\in[0,T]}|\f(t)|=\|\f(\cd)\|_\infty,\qq\forall\f(\cd)\in C([0,T];\dbR^n).\ea$$
Also, we define
$$\ba{ll}
\ns\ds L^{p+}(0,T;\dbR^n)=\bigcup_{q>p}L^q(0,T;\dbR^n),\qq p\in[1,\infty),\\
\ns\ds L^{p-}(0,T;\dbR^n)=\bigcap_{q<p}L^q(0,T;\dbR^n),\qq p\in(1,\infty).\ea$$
Next, for any continuous function $w:[0,T]\to[0,\infty)$, called a {\it weight function}, we define
$$\ba{ll}
\ns\ds L_{w(\cd)}^p(0,T;\dbR^n)=\Big\{\f:[0,T]\to\dbR^n\bigm|w(\cd)\f(\cd)\in L^p(0,T;\dbR^n)\Big\},\qq p\in[1,\infty],\\
\ns\ds C_{w(\cd)}([0,T];\dbR^n)=\Big\{\f:[0,T]\to\dbR^n\bigm|w(\cd)\f(\cd)\in C([0,T];\dbR^n)\Big\}.\ea$$
Clearly, if $\mbox{meas}\{t\in[0,T]\bigm|w(t)=0\}=0$, then $L^p_{w(\cd)}(0,T;\dbR^n)$ and $C_{w(\cd)}([0,T];\dbR^n)$ are normed linear spaces, under the following norms, respectively:
$$\|\f(\cd)\|_{L^p_{w(\cd)}}=\|w(\cd)\f(\cd)\|_p,\qq
\|\f(\cd)\|_{C_{w(\cd)}}=\|w(\cd)\f(\cd)\|_C.$$
Note that for any $\f(\cd)\in C_{w(\cd)}([0,T];\dbR^n)$, $\f(\cd)$ is continuous on the set $\{t\in[0,T]\bigm|w(t)>0\}$. If $w(s)=|s-s_0|^\g$ for some $s_0\in[0,T]$ and $\g>0$, then
$${a(\cd)\over|\cd-s_0|^\a}\in C_{w(\cd)}([0,T];\dbR^n),$$
for any $\a\les\g$, $a(\cd)\in C([0,T];\dbR^n)$. From the above, we should have some feeling about the space $C_{w(\cd)}([0,T];\dbR^n)$.

\ms

We denote
\bel{D}\D=\big\{(t,s)\in[0,T]^2\bigm|0\les s<t\les T\big\}.\ee
Note that the ``diagonal line'' $\{(t,t)\,|\,t\in[0,T]\}$ is not contained in $\D$. Thus if $\f:\D\to\dbR^n$ with $(t,s)\mapsto\f(t,s)$ being continuous, then $\f(\cd)$ could be unbounded as $|t-s|\to0$.

\ms

Before going further, let us first recall the Young's inequality for convolution (Theorem 3.9.4 in \cite{Bogachev 2007}).

\bl{Young} \sl Let $p,q,r\ges1$ satisfy
$${1\over p}+{1\over q}=1+{1\over r}.$$
Then for any $f(\cd)\in L^p(\dbR^n)$, $g(\cd)\in L^q(\dbR^n)$,
\bel{}\|f*g\|_r\les\|f\|_p\|g\|_q.\ee

\el

We now present several results. Some of them should be standard. However, we will provide the proofs for readers' convenience.

\bl{Lemma 2.1} \sl Let $\b\in(0,1)$ and $\f:\D\to\dbR^n$. Define
\bel{psi2.1}\psi(t)=\int_0^t{\f(t,s)\over(t-s)^{1-\b}}ds,\qq t\in[0,T].\ee

\ms

{\rm(i)} Suppose for some $p\in[1,\infty)$,
\bel{max|f|}\int_0^T\sup_{t\in[s,T]}|\f(t,s)|^pds<\infty.\ee
Then
\bel{|psi|}\|\psi(\cd)\|_p\les{T^\b\over\b}\(\int_0^T\sup_{t\in[s,T]}|\f(t,s)|^pds\)^{1\over p}.\ee

{\rm(ii)} Suppose, in addition, $t\mapsto\f(t,s)$ satisfies the following
\bel{}|\f(t,s)-\f(t',s)|\les\o(|t-t'|),\qq\forall t,t'\in[t_0-\si,t_0+\si],~s\in[0,t\land t'),\ee
for some modulus of continuity $\o:[0,\infty)\to[0,\infty)$, and for some $q>{1\over\b}$, $\si>0$, with $(t_0-\si,t_0+\si)\subseteq[0,T]$, the following holds:
\bel{int|f|}\int_{t_0-\si}^{t_0+\si}\sup_{t\in[s,t_0+\si]}|\f(t,s)|^qds<\infty.\ee
Then $\psi(\cd)$ is continuous at $t_0$. Consequently, if $t\mapsto\f(t,s)$ is continuous uniformly in $s\in[0,T]$ and \eqref{max|f|} holds for some $p>{1\over\b}$, then $t\mapsto\psi(t)$ is continuous on $[0,T]$.

\el

\it Proof. \rm (i) Let
$$\bar\f(s)=\sup_{t\in(s,T]}|\f(t,s)|,\qq s\in[0,T].$$
Then $\bar\f(\cd)\in L^p(0,T;\dbR)$. Hence, \eqref{|psi|} follows from Young's inequality for convolution.

\ms

(ii) Let $q>{1\over\b}$ which is equivalent to $\k\equiv(1-\b){q\over q-1}<1$ and let
$$\widehat{\f}(s)=\sup_{t\in(s,t_0+\si]}|\f(t,s)|,\qq s\in[0,t_0+\si].$$
For any $t_0-{\si\over m}<t<t'<t_0+{\si\over m}$ with $m\ges2$ large enough, we look at the following:
$$\ba{ll}
\ns\ds|\psi(t)-\psi(t')|=\Big|\int_0^t{\f(t,s)\over(t-s)^{1-\b}}ds-\int_0^{t'}{\f(t',s)\over(t'-s)^{1-\b}}ds\Big|
~~~~~~~~~~~~~~~~~~~~~~~~~~~~~~~~~~~~~~~~~~~~~~~~~\\
\ns\ds\les\int_0^{t-{\si\over m}}\Big|{\f(t,s)\over(t-s)^{1-\b}}-{\f(t',s)\over(t'-s)^{1-\b}}\Big|ds
+\int_{t-{\si\over m}}^t{|\f(t,s)|
\over(t-s)^{1-\b}}ds+\int_{t-{\si\over m}}^{t'}{|\f(t',s)|\over(t'-s)^{1-\b}}ds\\
\ns\ds\les\int_0^{t-{\si\over m}}|\f(t,s)|\({1\over(t-s)^{1-\b}}-{1\over(t'-s)^{1-\b}}\)ds
+\int_0^{t-{\si\over m}}{|\f(t,s)-\f(t',s)|\over(t'-s)^{1-\b}}ds\\
\ns\ds\qq+\int_{t-{\si\over m}}^t{\widehat{\f}(s)\over(t-s)^{1-\b}}ds+\int_{t-{\si\over m}}^{t'}{\widehat{\f}(s)\over(t'-s)^{1-\b}}ds\\
\ns\ds\les(t'-t)^{1-\b}\int_0^{t-{\si\over m}}{\widehat{\f}(s)\over(t-s)^{1-\b}(t'-s)^{1-\b}}ds
+\o(|t-t'|)\int_0^{t-{\si\over m}}{ds\over(t'-s)^{1-\b}}\\
\ns\ds\qq+\|\widehat{\f}(\cd)\|_{L^q(t-{\si\over m},t;\dbR^n)}\(\int_{t-{\si\over m}}^t{ds\over(t-s)^\k}\)^{q-1\over q}
+\|\widehat{\f}(\cd)\|_{L^q(t-{\si\over m},t';\dbR^n)}\(\int_{t-{\si\over m}}^{t'}{ds\over(t'-s)^\k}\)^{q-1\over q}\\
\ns\ds\les{(t'-t)^{1-\b}\over{({\si\over m})}^{2(1-\b)}}\|\widehat{\f}(\cd)\|_1
\1n+\o(|t\1n-\1n t'|){T^\b\over\b}\1n+\1n\|\widehat{\f}(\cd)\|_{L^q(t_0-{\si},t_0+{\si};\dbR^n)}
\[\({({\si\over m})^{1-\k}\over1-\k}\)^{q-1\over q}\3n+\1n\({(t'\1n-t+{\si\over m})^{1-\k}\over1-\k}\)^{q-1\over q}\].\ea$$
Hence, for any $\e>0$, we first take $m\ges1$ sufficiently large so that
$$\|\widehat{\f}(\cd)\|_{L^q(t_0-{\si},t_0+{\si};\dbR^n)}
\[\({{({\si\over m})}^{1-\k}\over1-\k}\)^{q-1\over q}+\({({3\si\over m})^{1-\k}\over1-\k}\)^{q-1\over q}\]<{\e\over2}.$$
Then let $\d\in(0,\si)$ be small enough so that
$${\d^{1-\b}\over{({\si\over m})}^{2(1-\b)}}\|\widehat{\f}(\cd)\|_1+\o(\d){T^\b\over\b}<{\e\over2}.$$
Combining the above, we see that $\psi(\cd)$ is continuous at $t_0$. The last conclusion follows easily from what we just proved. \endpf

\ms

The above lemma show that for any $p\in[1,\infty)$, under condition \eqref{max|f|}, one has $\psi(\cd)\in L^p(0,T;\dbR^n)$. To guarantee the continuity of $\psi(\cd)$ at $t_0\in(0,T]$, we need to assume the continuity of $t\mapsto\f(t,s)$ for $t\in[t_0-\sigma,t_0+\sigma]$, uniformly in $s\in[0,t)$, together with $L^q$ integrability of $\ds s\mapsto\sup_{t\in(s,t_0+\si]}|\f(t,s)|$ for $q>{1\over\b}$. The following example shows that continuity of $\psi(\cd)$ might fail $\bar\f(\cd)$ does not have a good enough integrability.

\bex{} \rm Let
$$\f(t,s)=\bar\f(s)={1\over|s-s_1|^{5\over6}},\q\b={1\over2},$$
with $s_1\in(0,T)$. Then $\bar\f(\cd)\in L^q(0,T;\dbR)$ with $q<{6\over5}<2={1\over\b}$. Note that
$$\psi(t)=\int_0^t{\f(t,s)\over(t-s)^{1-\b}}ds=\int_0^t{ds\over|s-s_1|^{5\over6}(t-s)^{1\over2}},\qq t\ges0.$$
By Young's inequality, we know that the above $\psi(\cd)\in L^p(0,T;\dbR)$ for some $p>1$. Also,
one sees that
$$\lim_{t\ua s_1}\psi(t)=\int_0^{s_1}{ds\over|s-s_1|^{4\over3}}=\infty.$$
Thus, $\psi(\cd)$ is discontinuous at $s_1$.

\ex

It is natural to ask if we relax the $L^q$-integrability of $\f(\cd)$, what can we say about the continuity of $\psi(\cd)$ defined by \eqref{psi2.1}? Let us make it more precise now. Let $\a_i,\b\in(0,1)$, $0\les i\les\ell$ and $0\les s_0<s_1<\cds<s_\ell\les T$. Define
\bel{w(s)}w(s)=\prod_{i=0}^\ell|s-s_i|^{1-\a_i},\qq s\in[0,T].\ee
Let $\wt\f:\D\to\dbR^n$ and define
\bel{wt psi}\wt\psi(t)=\int_0^t{\wt\f(t,s)\over w(s)(t-s)^{1-\b}}ds,\qq t\in[0,T].\ee
Comparing the above with \eqref{psi2.1}, we see that $\wt\psi(\cd)$ would be the same as $\psi(\cd)$ provided
\bel{f=f}\f(t,s)={\wt\f(t,s)\over w(s)},\qq(t,s)\in\D.\ee
We have the following result.

\bl{Lemma 2.2} \sl Let $\a_i,\b\in(0,1)$, $0\les i\les\ell$, and $w(\cd)$ be defined by \eqref{w(s)} with $0\les s_0<s_1<\cds<s_\ell\les T$. Let $\wt\f:\D\to\dbR^n$ satisfy
\bel{|f-f|}|\wt\f(t,s)-\wt\f(t',s)|\les\o(|t-t'|),\qq\forall(t,s),\ (t',s)\in\D,\ee
for some modulus of continuity $\o:[0,\infty)\to[0,\infty)$, and
\bel{int|f|<infty}\int_0^T\max_{t\in(s,T]}|\wt\f(t,s)|^qds<\infty,\ee
with some
\bel{q>}q>{1\over\b}\vee{1\over\a_i},\qq\q0\les i\les\ell.\ee
Define $\wt\psi(\cd)$ by \eqref{wt psi}. Then
\bel{}\wt\psi(\cd)\in L^\infty_{\bar w(\cd)}(0,T;\dbR^n)\bigcap\(\bigcap_{i=1}^\ell C\big((s_{i-1},s_i);\dbR^n\big)\),\ee
where
\bel{bar w}\bar w(s)=\prod_{i=0}^\ell|s-s_i|^{(1+{1\over q}-\a_i-\b)^+}.\ee
Consequently, for any $\e>0$,
\bel{2.16}\wt\psi(\cd)\in C_{\bar w^\e(\cd)}([0,T];\dbR^n),\ee
where
\bel{we(s)}\bar w^\e(s)=\prod_{i=0}^\ell|s-s_i|^{(1+{1\over q}-\a_i-\b)^++\e}.\ee
Further, if
\bel{a+b>1}\a_i+\b>1+{1\over q},\ee
then $\wt\psi(\cd)$ is continuous at $s_i$. Consequently, if
\bel{a+b>1*}\min_{0\les i\les\ell}\a_i+\b>1+{1\over q},\ee
then $\wt\psi(\cd)\in C([0,T];\dbR^n)$.

\el

\it Proof. \rm Denote
$$\bar\f(s)=\max_{t\in[s,T]}|\wt\f(t,s)|,\qq s\in[0,T].$$
Define
$$w_i(s)=\prod_{j\ne i}|s-s_j|^{1-\a_j},\qq s\in[0,T],~0\les i\les\ell,$$
and
\bel{d0}\d_0=\min_{0\les i\les\ell-1}{s_{i+1}-s_i\over2}>0,\qq\bar s_i={s_i+s_{i+1}\over2},\q0\les i\les\ell-1.\ee
Then
\bel{s+d<s<s-d}s_i+\d_0\les\bar s_i\les s_{i+1}-\d_0,\qq i=0,1,\cds,\ell-1.\ee
We establish some estimates for $\widetilde{\psi}(\cd)$ on some time intervals, more precisely, on $[0,s_0)$, $(s_0,\bar s_0]$,
$(0,\bar s_0]$, $[\bar s_0,s_1)$, $(s_1,\bar s_1]$, $[\bar s_1,s_2)$. Then induction will apply.

\ms
(i) For $t\in[0,s_0)$,  one has
$$\ba{ll}
\ns\ds|\wt\psi(t)|=\Big|\int_0^t{\wt\f(t,s)\over w_0(s)(s_0-s)^{1-\a_0}(t-s)^{1-\b}}ds\Big|\les K\int_0^t{\bar\f(s)\over (s_0-s)^{1-\a_0}(t-s)^{1-\b}}ds\\
\ns\ds\qq~\les K\(\int_0^t\bar\f(s)^qds\)^{1\over q}\(\int_0^t{ds\over (s_0-s)^{(1-\a_0){q\over q-1}}(t-s)^{(1-\b){q\over q-1}}}\)^{q-1\over q}\\
\ns\ds\qq~\les K\(\int_0^t\bar\f(s)^qds\)^{1\over q}(s_0-t)^{\a_0+\b-2}\(\int_0^t{ds\over(1+{t-s\over s_0-t})^{(1-\a_0){q\over q-1}}({t-s\over s_0-t})^{(1-\b){q\over q-1}}}\)^{q-1\over q}.\ea$$
Here and throughout the paper, $K$ is a positive constant, which may be different when appears at different places.
Now, we let $\t={t-s\over s_0-t}$. Then $s=t-(s_0-t)\t$, $ds=(t-s_0)d\t$, and
$$\ba{ll}
\ns\ds|\wt\psi(t)|\les K\(\int_0^t\bar\f(s)^qds\)^{1\over q}(s_0-t)^{\a_0+\b-1-{1\over q}}\(\int_0^{t\over s_0-t}{d\t\over(1+\t)^{(1-\a_0){q\over q-1}}
\t^{(1-\b){q\over q-1}}}\)^{q-1\over q}\\
\ns\ds\qq~\leq K(s_0-t)^{\a_0+\b-1-{1\over q}}\[\(\int_0^1{d\t\over\t^{(1-\b){q\over q-1}}}\)^{q-1\over q}
+\(\int_1^{s_0\over s_0-t}{d\t\over\t^{(2-\a_0-\b){q\over q-1}}}\)^{q-1\over q}\]\\%
\ns\ds\qq~\les K(s_0-t)^{\a_0+\b-1-{1\over q}}\[1+\({s_0\over s_0-t}\)^{[1-(2-\a_0-\b){q\over q-1}]{q-1\over q}}\]\\
\ns\ds\qq~\les K(s_0-t)^{\a_0+\b-1-{1\over q}}\[1+\({1\over s_0-t}\)^{\a_0+\b-1-{1\over q}}\]\les K+K(s_0-t)^{\a_0+\b-1-{1\over q}},\q t\in[0,s_0).\ea$$
Hence,
\bel{}(s_0-t)^{(1-\a_0-\b+{1\over q})^+}|\wt\psi(t)|\les K,\qq\q t\in[0, s_0).\ee

(ii) For $t\in(s_0,\bar s_0]$, one has
$$\ba{ll}
\ns\ds|\wt\psi(t)|\les\int_0^{s_0}{\bar\f(s)\over w(s)(t-s)^{1-\b}}ds+\int_{s_0}^t{\bar\f(s)\over w(s)(t-s)^{1-\b}}ds\equiv I_1+I_2.\ea$$
For $I_1$, we have
$$\ba{ll}
\ns\ds I_1~\les K\(\int_0^{s_0}{ds\over (s_0-s)^{(1-\a_0){q\over q-1}}(t-s)^{(1-\b){q\over q-1}}}\)^{q-1\over q}\\
\ns\ds\q\les K(t-s_0)^{\a_0+\b-2}\(\int_0^{s_0}{ds\over({s_0-s\over t-s_0})^{(1-\a_0){q\over q-1}}
(1+{s_0-s\over t-s_0})^{(1-\b){q\over q-1}}}\)^{q-1\over q}.\ea$$
Now, we let $\t={s_0-s\over t-s_0}$. Then $s=s_0-(t-s_0)\t$, $ds=(s_0-t)d\t$, and
$$\ba{ll}
\ns\ds I_1~\les K(t-s_0)^{\a_0+\b-1-{1\over q}}\(\int_0^{s_0\over t-s_0}{d\t\over\t^{{(1-\a_0){q\over q-1}}}(1+\t)^{(1-\b){q\over q-1}}}\)^{q-1\over q}\\
\ns\ds\qq\les K(t-s_0)^{\a_0+\b-1-{1\over q}}\[\(\int_0^1{d\t\over\t^{(1-\a_0){q\over q-1}}}\)^{q-1\over q}
+\(\int_1^{s_0\over t-s_0}{d\t\over\t^{(2-\a_0-\b){q\over q-1}}}\)^{q-1\over q}\]\\
\ns\ds\q\les K(t-s_0)^{\a_0+\b-1-{1\over q}}\[1+\({s_0\over t-s_0}\)^{[1-(2-\a_0-\b){q\over q-1}]{q-1\over q}}\]\\
\ns\ds\q\les K(t-s_0)^{\a_0+\b-1-{1\over q}}\[1+\({1\over t-s_0}\)^{\a_0+\b-1-{1\over q}}\]\les K+K(t-s_0)^{\a_0+\b-1-{1\over q}},\q t\in(s_0,\bar s_0].\ea$$
Now, we look at $I_2$, noting $s_0+\d_0\les\bar s_0\les s_1-\d_0$,
$$\ba{ll}
\ns\ds I_2~=\int_{s_0}^t{\bar\f(s)\over w(s)(t-s)^{1-\b}}ds=\int_{s_0}^t{\bar\f(s)\over w_0(s)(s-s_0)^{1-\a_0}(t-s)^{1-\b}}ds\\
\ns\ds\q\les K\(\int_{s_0}^t\bar\f(s)^qds\)^{1\over q}\(\int_{s_0}^t{ds\over (s-s_0)^{(1-\a_0){q\over q-1}}(t-s)^{(1-\b){q\over q-1}}}\)^{q-1\over q}\\
\ns\ds\q\les K(t-s_0)^{\a_0+\b-2}
\(\int_{s_0}^t{ds\over({s-s_0\over t-s_0})^{(1-\a_0){q\over q-1}}({1-{s-s_0\over t-s_0}})^{(1-\b){q\over q-1}}}\)^{q-1\over q}.\ea$$
 Note that by \eqref{q>}, we have
$$(1-\a_0){q\over q-1}<1,\q(1-\b){q\over q-1}<1,$$
which are equivalent to the following:
$$1-(1-\a_0){q\over q-1}={\a_0q-1\over q-1}>0,\qq
1-(1-\b){q\over q-1}={\b q-1\over q-1}>0.$$
Let $\t={s-s_0\over t-s_0}$. Then $s=s_0+(t-s_0)\t$, $ds=(t-s_0)d\t$, and
$$\ba{ll}
\ns\ds\qq I_2\les K(t-s_0)^{\a_0+\b-1-{1\over q}}\(\int_0^1
{d\t\over\t^{1-{\a_0q-1\over q-1}}(1-\t)^{1-{\b q-1\over q-1}}}\)^{q-1\over q}=K(t-s_0)^{\a_0+\b-1-{1\over q}}B\({\a_0q-1\over q-1},{\b q-1\over q-1}\)^{q-1\over q}.\ea$$
Here, $(a,b)\mapsto B(a,b)$ is the Beta function. Hence,
\bel{}(t-s_0)^{(1-\a_0-\b+{1\over q})^+}|\wt\psi(t)|\les K,\qq\q t\in(s_0,\bar s_0].\ee

(iii) For $t\in[\bar s_0,s_1)$, we have
$$\ba{ll}
\ns\ds|\wt\psi(t)|\les\Big|\int_0^{\bar s_0}{\wt\f(t,s)\over w(s)(t-s)^{1-\b}}ds\Big|
+\Big|\int_{\bar s_0}^t{\wt\f(t,s)\over w(s)(t-s)^{1-\b}}ds\Big|\equiv I_3+I_4.\ea$$
For $I_3$, we have
$$\ba{ll}
\ns\ds I_3\les\int_{0}^{s_0}{\bar\f(s)\over w(s)(t-s)^{1-\b}}ds+\int_{s_0}^{\bar s_0}{\bar\f(s)\over w(s)(t-s)^{1-\b}}ds\\
\ns\ds\q\les\int_{0}^{s_0}{\bar\f(s)\over (s-s_0)^{1-\a_0}}ds+\int_{s_0}^{\bar s_0}{\bar\f(s)\over (s-s_0)^{1-\a_0}(\bar s_0-s)^{1-\b}}ds\\
\ns\ds\q\les K\(\int_0^{s_0}\bar\f(s)^qds\)^{1\over q}\(\int_0^{s_0}{ {ds\over(s-s_0)^{(1-\a_0){q\over q-1}}}\)^{q-1\over q}}\\
\ns\ds\qq+{{K\(\int_{s_0}^{\bar s_0}\bar\f(s)^qds\)^{1\over q}\(\int_{s_0}^{\bar s_0}{ds\over (s-s_0)^{(1-\a_0){q\over q-1}}(\bar s_0-s)^{(1-\b){q\over q-1}}}}\)^{q-1\over q}}\\
\ns\ds\q\les K+ K(\bar s_0-s_0)^{\a_0+\b-2}
\(\int_{s_0}^{\bar s_0}{ds\over({s-s_0\over \bar s_0-s_0})^{(1-\a_0){q\over q-1}}({1-{s-s_0\over \bar s_0-s_0}})^{(1-\b){q\over q-1}}}\)^{q-1\over q}.\ea$$
Let $\t={s-s_0\over \bar s_0-s_0}$. Then $s=s_0+(\bar s_0-s_0)\t$, $ds=(\bar s_0-s_0)d\t$, and
$$\ba{ll}
\ns\ds I_3\les K+K(\bar s_0-s_0)^{\a_0+\b-1-{1\over q}}\(\int_0^1
{d\t\over\t^{1-{\a_0q-1\over q-1}}(1-\t)^{1-{\b q-1\over q-1}}}\)^{q-1\over q}\\
\ns\ds\q=K+K(\bar s_0-s_0)^{\a_0+\b-1-{1\over q}}B\({\a_0q-1\over q-1},{\b q-1\over q-1}\)^{q-1\over q}\les K.\ea$$
For $I_4$, we have
$$\ba{ll}
\ns\ds I_4
\les\int_{\bar s_0}^t{\bar\f(s)\over w_1(s)(s_1-s)^{1-\a_1}(t-s)^{1-\b}}ds\\
\ns\ds\q\les K\(\int_{\bar s_0}^t\bar\f(s)^qds\)^{1\over q}\(\int_{\bar s_0}^t{ds\over (s_1-s)^{(1-\a_1){q\over q-1}}(t-s)^{(1-\b){q\over q-1}}}ds\)^{q-1\over q}\\
\ns\ds\q\les K(s_1-t)^{\a_1+\b-2}
\(\int_{\bar s_0}^t{ds\over(1+{t-s\over s_1-t})^{(1-\a_1){q\over q-1}}({t-s\over s_1-t})^{(1-\b){q\over q-1}}}ds\)^{q-1\over q}.\ea$$
Let $\t={t-s\over s_1-t}$. Then $s=t-(s_1-t)\t$, $ds=-(s_1-t)d\t$, and
$$\ba{ll}
\ns\ds I_4\les K(s_1-t)^{\a_1+\b-1-{1\over q}}\(\int_0^{t-\bar s_0\over s_1-t}{d\t\over(1+\t)^{(1-\a_1){q\over q-1}}\t^{(1-\b){q\over q-1}}}\)^{q-1\over q}\\
\ns\ds\q\les K(s_1-t)^{\a_1+\b-1-{1\over q}}\[\(\int_0^1{d\t\over\t^{(1-\b){q\over q-1}}}\)^{q-1\over q}
+\(\int_1^{s_1-\bar s_0\over s_1-t}{d\t\over\t^{(2-\a_1-\b){q\over q-1}}}\)^{q-1\over q}\]\\
\ns\ds\q\les K(s_1-t)^{\a_1+\b-1-{1\over q}}\[1+\({s_1-\bar s_0\over s_1-t}\)^{[1-(2-\a_1-\b){q\over q-1}]{q-1\over q}}\]\\
\ns\ds\q\les K(s_1-t)^{\a_1+\b-1-{1\over q}}\[1+\({1\over s_1-t}\)^{\a_1+\b-1-{1\over q}}\]\les K(s_1-t)^{\a_1+\b-1-{1\over q}},\qq t\in[\bar s_0,s_1).\ea$$
Then we have
\bel{}(s_1-t)^{(1-\a_1-\b+{1\over q})^+}|\wt\psi(t)|\les K,\qq t\in[\bar s_0,s_1).\ee

(iv) For $t\in(s_1,\bar s_1]$, we have
$$\ba{ll}
\ns\ds|\wt\psi(t)|\les\int_0^{s_1}{\bar\f(s)\over w(s)(t-s)^{1-\b}}ds+\int_{s_1}^t{\bar\f(s)\over w(s)(t-s)^{1-\b}}ds\equiv I_5+I_6.\ea$$
Note that $\bar s_0\les s_1-\d_0$,
$$\ba{ll}
\ns\ds I_5=\1n\int_0^{s_1}\2n{\bar\f(s)\over w(s)(t-s)^{1-\b}}ds\1n=\1n\int_0^{\bar s_0}\2n{\bar\f(s)\over w_0(s)|s-s_0|^{1-\a_0}(t-s)^{1-\b}}ds\1n+\1n\int_{\bar s_0}^{s_1}\2n{\bar\f(s)\over w_1(s)(s_1-s)^{1-\a_1}(t-s)^{1-\b}}ds\\
\ns\ds\q\les K\int_0^{\bar s_0}{\bar\f(s)\over |s-s_0|^{1-\a_0}}ds+K\int_{\bar s_0}^{s_1}{\bar\f(s)\over(s_1-s)^{1-\a_1}[(t-s_1)+(s_1-s)]^{1-\b}}ds\\
\ns\ds\q\les K\(\int_0^{\bar s_0}\bar\f(s)^qds\)^{1\over q}\[\(\int_0^{\bar s_0}{ds\over |s-s_0|^{(1-\a_0){q\over q-1}}}\)^{q-1\over q}\\
\ns\ds\qq\qq+(t-s_1)^{\a_1+\b-2}\(\int_{\bar s_0}^{s_1}{ds\over({s_1-s\over t-s_1})^{(1-\a_1){q\over q-1}}(1+{s_1-s\over t-s_1})^{(1-\b){q\over q-1}}}\)^{q-1\over q}\].\ea$$
Let $\t={s_1-s\over t-s_1}$. Then $s=s_1-(t-s_1)\t$, $ds=-(t-s_1)d\t$, and
$$\ba{ll}
\ns\ds I_5\les K\[1+(t-s_1)^{\a_1+\b-2}\(\int_{\bar s_0}^{s_1}{ds\over({s_1-s\over t-s_1})^{(1-\a_1){q\over q-1}}(1+{s_1-s\over t-s_1})^{(1-\b){q\over q-1}}}\)^{q-1\over q}\]\\
\ns\ds\q\les K\[1+(t-s_1)^{\a_1+\b-1-{1\over q}}\(\int_0^{s_1-\bar s_0\over t-s_1}{d\t\over\t^{(1-\a_1){q\over q-1}}(1+\t)^{(1-\b){q\over q-1}}}\)^{q-1\over q}\]\\
\ns\ds\q\les K\[1+(t-s_1)^{\a_1+\b-1-{1\over q}}\(\int_0^1{d\t\over\t^{(1-\a_1){q\over q-1}}}+\int_1^{s_1-\bar s_0\over t-s_1}{d\t\over\t^{(2-\a_1-\b){q\over q-1}}}\)^{q-1\over q}\]\\
\ns\ds\q\les K\[1+(t-s_1)^{\a_1+\b-1-{1\over q}}+(t-s_1)^{\a_1+\b-1-{1\over q}}\({s_1-\bar s_0\over t-s_1}\)^{[1-(2-\a_1-\b){q\over q-1}]{q-1\over q}}\]\\
\ns\ds\q\les K+K(t-s_1)^{\a_1+\b-1-{1\over q}}.\ea$$
Now, we look at $I_6$, noting $s_1+\d_0\les\bar s_1\les s_2-\d_0$,
$$\ba{ll}
\ns\ds I_6=\int_{s_1}^t{\bar\f(s)\over w(s)(t-s)^{1-\b}}ds=\int_{s_1}^t{\bar\f(s)\over w_1(s)(s-s_1)^{1-\a_1}(t-s)^{1-\b}}ds\\
\ns\ds\q\les K\(\int_{s_1}^t\bar\f(s)^qds\)^{1\over q}\(\int_{s_1}^t{ds\over(s-s_1)^{(1-\a_1){q\over q-1}}(t-s)^{(1-\b){q\over q-1}}}\)^{q-1\over q}\\
\ns\ds\q\les K(t-s_1)^{\a_1+\b-2}\(\int_{s_1}^t{ds\over({s-s_1\over t-s_1})^{(1-\a_1){q\over q-1}}(1-{s-s_1\over t-s_1})^{(1-\b){q\over q-1}}}\)^{q-1\over q}.\ea$$
Let $\t={s-s_1\over t-s_1}$. Then $s=s_1+(t-s_1)\t$, $ds=(t-s_1)d\t$, and
$$\ba{ll}
\ns\ds I_6\1n\les\1n K(t\1n-\1n s_1)^{\a_1+\b-1-{1\over q}}\(\int_0^1\2n{d\t\over\t^{1-{\a_1q-1\over q-1}}(1-\t)^{1-{\b q-1\over q-1}}}\)^{q-1\over q}\3n=\1n K(t-s_1)^{\a_1+\b-1-{1\over q}}B\({\a_1q-1\over q-1},{\b q-1\over q-1}\)^{q-1\over q}.\ea$$
Hence,
\bel{}(t-s_1)^{(1+{1\over q}-\a_1-\b)^+}|\wt\psi(t)|\les K,\qq t\in(s_1,\bar s_1].\ee

(v) For $t\in[\bar s_1,s_2)$,
$$|\wt\psi(t)|\les\int_0^{s_1}{\bar\f(s)\over w(s)(t-s)^{1-\b}}ds+\int_{s_1}^{\bar s_1}{\bar\f(s)\over w(s)(t-s)^{1-\b}}ds
+\int_{\bar s_1}^t{\bar\f(s)\over w(s)(t-s)^{1-\b}}ds\equiv I_7+I_8+I_9.$$
We look at the three terms one-by-one. Since $\bar s_1\ges s_1+\d_0$, one has
$$\ba{ll}
\ns\ds I_7=\int_0^{s_1}{\bar\f(s)\over w(s)(t-s)^{1-\b}}ds\les K\Big(\int_0^{\bar s_0}{\bar\f(s)\over |s-s_0|^{1-\a_0}(t-s)^{1-\b}}ds+\int_{\bar s_0}^{s_1}
{\bar\f(s)\over (s_1-s)^{1-\a_1}(t-s)^{1-\b}}ds\Big)\\
\ns\ds\q\les\1n K\(\int_0^{\bar s_0}\2n\bar\f(s)^qds\)^{1\over q}\(\int_0^{\bar s_0}\2n{ds\over |s-s_0|^{(1-\a_0){q\over q-1}}}\Big)^{q-1\over q}\3n+\1n K\(\int_{\bar s_0}^{s_1}\bar\f(s)^qds\)^{1\over q}\(\int_{\bar s_0}^{s_1}\2n{ds\over (s_1-s)^{(1-\a_1){q\over q-1}}}\)^{q-1\over q}\3n\les\1n K.\ea$$
For $I_8$, since $(t-s)\ges(s_1-s)$, one has
$$\ba{ll}
\ns\ds I_8\1n=\2n\int_{s_1}^{\bar s_1}\2n{\bar\f(s)\over w_1(s)(s-s_1)^{1-\a_1}(t-s)^{1-\b}}ds\1n\les\1n K\(\int_{s_1}^{\bar s_1}\2n\bar\f(s)^qds\)^{1\over q}\(\int_{s_1}^{\bar s_1}\2n{ds\over(s-s_1)^{(1-\a_1){q\over q-1}}(\bar s_1-s)^{(1-\b){q\over q-1}}}\)^{q-1\over q}\\
\ns\ds\q\les K(\bar s_1-s_1)^{\a_1+\b-2}\(\int_{s_1}^{\bar s_1}{ds\over({s-s_1\over\bar s_1-s_1})^{(1-\a_1){q\over q-1}}({\bar s_1-s\over\bar s_1-s_1})^{(1-\b){q\over q-1}}}\)^{q-1\over q}.\ea$$
Let $\t={s-s_1\over\bar s_1-s_1}$. Then $s=s_1+(\bar s_1-s_1)\t$, $ds=(\bar s_1-s_1)d\t$, and
$$\ba{ll}
\ns\ds I_8\les K(\bar s_1-s_1)^{\a_1+\b-1-{1\over q}}\(\int_0^1{d\t\over\t^{1-{\a_1q-1\over q-1}}(1-\t)^{1-{\b q-1\over q-1}}}\)^{q-1\over q}\les K.\ea$$
Finally, for $I_9$, one has
$$\ba{ll}
\ns\ds I_9\1n=\2n\int_{\bar s_1}^t\2n{\bar\f(s)\over w_2(s)(s_2-s)^{1-\a_2}(t-s)^{1-\b}}ds\1n\les\1n K\(\int_{\bar s_1}^t\2n\bar\f(s)^qds\)^{1\over q}\(\int_{\bar s_1}^t\2n{ds\over (s_2-s)^{(1-\a_2){q\over q-1}}(t-s)^{(1-\b){q\over q-1}}}\)^{q-1\over q}\\
\ns\ds\q\les K(s_2-t)^{\a_2+\b-2}
\(\int_{\bar s_1}^t{ds\over(1+{t-s\over s_2-t})^{(1-\a_2){q\over q-1}}({t-s\over s_2-t})^{(1-\b){q\over q-1}}}\)^{q-1\over q}.\ea$$
Let $\t={t-s\over s_2-t}$. Then $s=t-(s_2-t)\t$, $ds=-(s_2-t)d\t$, and
$$\ba{ll}
\ns\ds I_9=\int_{\bar s_1}^t{\bar\f(s)\over w(s)(t-s)^{1-\b}}ds\les K(s_2-t)^{\a_2+\b-1-{1\over q}}\(\int_0^{t-\bar s_1\over s_2-t}{d\t\over(1+\t)^{(1-\a_2){q\over q-1}}\t^{(1-\b){q\over q-1}}}\)^{q-1\over q}\\
\ns\ds\q\les K(s_2-t)^{\a_2+\b-1-{1\over q}}\(\int_0^1{ds\over\t^{(1-\b){q\over q-1}}}ds+
\int_1^{s_2-\bar s_1\over s_2-t}{ds\over\t^{(2-\a_2-\b){q\over q-1}}}\)^{q-1\over q}\\
\ns\ds\q\les K(s_2-t)^{\a_2+\b-1-{1\over q}}\[1+\({s_2-\bar s_1\over s_2-t}\)^{[1-(2-\a_2-\b){q\over q-1}]{q-1\over q}}\]\les K+K(s_2-t)^{\a_2+\b-1-{1\over q}}.\ea$$
Hence,
\bel{}(s_2-t)^{(1-\a_2-\b+{1\over q})^+}|\wt\psi(t)|\les K,\qq t\in[\bar s_1,s_2).\ee
By induction, we can obtain the following:
\bel{}\left\{\2n\ba{ll}
\ns\ds(s_0-t)^{(1+{1\over q}-\a_0-\b)^+}|\wt\psi(t)|\les K,\qq t\in(0,s_0),\\
\ns\ds(t-s_i)^{(1+{1\over q}-\a_i-\b)^+}|\wt\psi(t)|\les K,\qq t\in(s_i,\bar s_i],\qq0\les i\les\ell-1,\\
\ns\ds(s_{i+1}-t)^{(1+{1\over q}-\a_{i+1}-\b)^+}|\wt\psi(t)|\les K,\qq t\in[\bar s_i,s_{i+1}),\qq0\les i\les\ell-1,\\
\ns\ds(t-s_l)^{(1+{1\over q}-\a_l-\b)^+}|\wt\psi(t)|\les K,\qq t\in(s_l,T].\ea\right.\ee
Hence, for $\bar w(\cd)$, we have
$$\bar w(t)|\wt\psi(t)|\les K,\qq\forall t\in[0,T]\setminus\{s_0,s_1,\cds,s_\ell\}.$$
Then \eqref{2.16} follows easily. On the other hand, if $s_0=0$ and $s_l=T$,
then  for any $t_0\in[0,T]\setminus\{s_0,s_1,s_2,\cds,s_\ell\}$, one can find a $\si>0$  such that for some $i=0,1,2,\cds,\ell-1$, $[t_0-\si,t_0+\si]\subseteq(s_i,s_{i+1})$. Then
$$\int_{t_0-\si}^{t_0+\si}\sup_{t\in[s,t_0+\si]}\Big|{\wt\f(t,s)\over w(s)}\Big|^qds
\les K\int_{t_0-\si}^{t_0+\si}\bar\f(s)^qds<\infty.$$
Hence, by the similar argument of Lemma \ref{Lemma 2.1}, $\wt\psi(\cd)$ is continuous at any such a $t_0$. If $0<s_0$ or $s_\ell<T$, we can use the similar argument to show the continuity of $\wt\psi(\cd)$ at $t_0\in[0,T]\setminus\{s_0,s_1,s_2,\cds,s_\ell\}$.

\ms

Finally, if \eqref{a+b>1} holds, then for $\si>0$ small enough, and for $r<q$, with $1+{1\over q}-\a_i<{1\over r}<\b$,
$$\ba{ll}
\ns\ds\int_{s_i-\si}^{s_i+\si}{\bar\f(s)^r\over w(s)^r}ds\les K\(\int_{s_i-\si}^{s_i+\si}\bar\f(s)^qds\)^{r\over q}\(\int_{s_i-\si}^{s_i+\si}{ds\over|s-s_i|^{(1-\a_i){rq\over q-r}}}\)^{q-r\over q}<\infty,\ea$$
since
$$\q1+{1\over q}-\a_i<{1\over r}\q\iff\q(1-\a_i){rq\over q-r}<1.$$
Thus,  by the similar argument of Lemma \ref{Lemma 2.1}, we have the continuity of $\wt\psi(\cd)$ at $s_i$. The last conclusion is clear. \endpf

\ms

From the above, we see that if $\a_i+\b>1$ for all $0\les i\les\ell$ and
\bel{}\bar\f(\cd)\in L^{{1\over\a_i+\b-1}+}(0,T;\dbR),\qq0\les i\les\ell,\ee
then $\wt\psi(\cd)\in C([0,T];\dbR^n)$. On the other hand, if $\a_i+\b<1$, and $\bar\f(\cd)$ is essentially non-zero near $s_i$, then $\wt\psi(\cd)$ will be blow-up near $s_i$, and roughly it will grow no more than $|t-s_i|^{\a_i+\b-1-{1\over q}}$.

\ms

The following result is a kind of Gronwall's inequality with a singular kernel.


\bl{Gronwall} \sl Let $\b\in(0,1)$ and $q>{1\over\b}$. Let $L(\cd),a(\cd),y(\cd)$ be nonnegative functions with
$$L(\cd)\in L^q(0,T;\dbR),\q a(\cd),\ y(\cd)\in L^{q\over q-1}(0,T;\dbR).$$
Suppose
\bel{2.2}y(t)\les a(t)+\int_0^t{L(s)y(s)\over(t-s)^{1-\b}}ds,\qq \ae~t\in[0,T].\ee
Then
\bel{y<}y(t)\les a(t)+\sum_{i=0}^{k-1}c_i\int_0^t{L(s)a(s)\over(t-s)^{1-\b_i}}ds
+c_k\int_0^tL(s)a(s)ds,\qq\ae~t\in[0,T],\ee
for some constants $c_i>0$ and $\b_i\in(0,1)$ defined by
$$\b_i=\b+i\(\b-{1\over q}\),\qq0\les i\les k-1,$$
with $k$ being the smallest integer satisfying
$$\b+k\(\b-{1\over q}\)\ges1.$$

\el

\it Proof. \rm First of all, since $L(\cd)\in L^q(0,T;\dbR)$ and $y(\cd)\in L^{q'}(0,T;\dbR)$, $q'={q\over q-1}$, we know that $L(\cd)y(\cd)\in L^1(0,T;\dbR)$. Hence, the integral on the right hand side of \eqref{2.2} is well-defined, as a function in $L^1(0,T;\dbR)$. Now, we observe the following:
$$\ba{ll}
\ns\ds y(t)\les a(t)+\int_0^t{L(s)y(s)\over(t-s)^{1-\b}}ds\les a(t)+\int_0^t{L(s)a(s)\over(t-s)^{1-\b}}ds+\int_0^t{L(s)\over (t-s)^{1-\b}}\int_0^s{L(\t)y(\t)\over(s-\t)^{1-\b}}d\t ds\\
\ns\ds\qq\les a(t)+\int_0^t{L(s)a(s)\over (t-s)^{1-\b}}ds+\int_0^tL(\t)\[\int_\t^t{L(s)\over (t-s)^{1-\b}(s-\t)^{1-\b}}ds\]y(\t)d\t.\ea$$
Let $\ds r={s-\t\over t-\t}$. Then $s=\t+(t-\t)r$ and $ds=(t-\t)dr$. Thus,
$$\ba{ll}
\ns\ds\int_\t^t{L(s)\over(t-s)^{1-\b}(s-\t)^{1-\b}}ds\les\(\int_\t^tL(s)^qds\)^{1\over q}\(\int_\t^t{ds\over(t-s)^{(1-\b)q'}(s-\t)^{(1-\b)q'}}ds\)^{1\over q'}\\
\ns\ds\les\|L(\cd)\|_q\(\int_0^1{(t-\t)dr\over[(t-\t)(1-r)]^{(1-\b)q'}[(t-\t)r]^{(1-\b)q'}}\)^{1\over q'}\\
\ns\ds=\|L(\cd)\|_q{1\over(t-\t)^{2(1-\b)-{1\over q'}}}\(\int_0^1{dr\over(1-r)^{(1-\b)q'}r^{(1-\b)q'}}\)^{1\over q'}.\ea$$
Since $q>{1\over\b}$ which is equivalent to
\bel{<1}0<1-(1-\b)q'=1-{(1-\b)q\over q-1}={q-1-q+\b q\over q-1}={\b q-1\over q-1},\ee
 we obtain
$$\ba{ll}
\ns\ds\int_\t^t{L(s)\over(t-s)^{1-\b}(s-\t)^{1-\b}}ds\les{\|L(\cd)\|_q\over
(t-\t)^{2(1-\b)-{1\over q'}}}B\({\b q-1\over q-1},{\b q-1\over q-1}\)^{1\over q'}\equiv{c_1\over(t-\t)^{1-\b_1}},\ea$$
with $B(\cd\,,\cd)$ being the Beta function and
$$\ba{ll}
\ns\ds c_1=\|L(\cd)\|_qB\({\b q-1\over q-1},{\b q-1\over q-1}\)^{1\over q'},\\
\ns\ds\b_1=1-\(2(1-\b)-{1\over q'}\)=2\b+{1\over q'}-1=\b+\(\b-{1\over q}\)>\b.\ea$$
Consequently,
$$\ba{ll}
\ns\ds y(t)\les a(t)+\int_0^t{L(s)a(s)\over(t-s)^{1-\b}}ds+c_1\int_0^t{L(s)y(s)\over(t-s)^{1-\b_1}}ds\\
\ns\ds\qq\les a(t)+\int_0^t{L(s)a(s)\over(t-s)^{1-\b}}ds+c_1\int_0^t{L(s)a(s)\over(t-s)^{1-\b_1}}ds+c_1\int_0^t{L(s)\over(t-s)^{1-\b_1}}\int_0^s{L(\t)y(\t)\over(s-\t)^{1-\b}}d\t ds\\
\ns\ds\qq=a(t)+\sum_{i=0}^1c_i\int_0^t{L(s)a(s)\over(t-s)^{1-\b_i}}ds+c_1\int_0^tL(\t)\[\int_\t^t{L(s)\over
(t-s)^{1-\b_1}(s-\t)^{1-\b}}ds\]y(\t)d\t,\ea$$
with $c_0=1$ and $\b_0=\b$. Let $\ds r={s-\t\over t-\t}$. Then $s=\t+(t-\t)r$ and $ds=(t-\t)dr$. Thus,
$$\ba{ll}
\ns\ds\int_\t^t{L(s)\over(t-s)^{1-\b_1}(s-\t)^{1-\b}}ds\les\(\int_\t^tL(s)^qds\)^{1\over q}\(\int_\t^t{ds\over(t-s)^{(1-\b_1)q'}(s-\t)^{(1-\b)q'}}ds\)^{1\over q'}\\
\ns\ds\les\|L(\cd)\|_q\(\int_0^1{(t-\t)dr\over[(t-\t)(1-r)]^{(1-\b_1)q'}[(t-\t)r]^{(1-\b)q'}}\)^{1\over q'}\\
\ns\ds=\|L(\cd)\|_q{1\over(t-\t)^{2-\b-\b_1-{1\over q'}}}\(\int_0^1{dr\over(1-r)^{(1-\b_1)q'}r^{(1-\b)q'}}\)^{1\over q'}.\ea$$
Since $q>{1\over\b}>{1\over\b_1}$, we have
\bel{<1*}0<1-(1-\b_1)q'=1-{(1-\b_1)q\over q-1}={q-1-q+\b_1q\over q-1}={\b_1q-1\over q-1}.\ee
Hence, we obtain
$$\ba{ll}
\ns\ds c_1\int_\t^t{L(s)\over(t-s)^{1-\b_1}(s-\t)^{1-\b}}ds\les{c_1\|L(\cd)\|_q\over
(t-\t)^{2-\b-\b_1-{1\over q'}}}B\({\b q-1\over q-1},{\b_1 q-1\over q-1}\)^{1\over q'}\equiv{c_2\over(t-\t)^{1-\b_2}},\ea$$
with
$$\ba{ll}
\ns\ds c_2=c_1\|L(\cd)\|_qB\({\b q-1\over q-1},{\b_1 q-1\over q-1}\)^{1\over q'},\\
\ns\ds\b_2=1-\(2-\b-\b_1-{1\over q'}\)=\b+\b_1+{1\over q'}-1=\b+2\(\b-{1\over q}\)>\b.\ea$$
Consequently,
$$y(t)\les a(t)+\sum_{i=0}^1c_i\int_0^t{L(s)a(s)\over(t-s)^{1-\b_i}}ds+c_2\int_0^t{L(s)y(s)\over
(t-s)^{1-\b_2}}ds.$$
By induction, we are able to show that
$$y(t)\les a(t)+\sum_{i=0}^{k-1}c_i\int_0^t{L(s)a(s)\over(t-s)^{1-\b_i}}ds
+c_k\int_0^t{L(s)y(s)\over(t-s)^{1-\b_k}}ds,$$
with
$$\b_i=\b+i\(\b-{1\over q}\),\qq0\les i\les k,$$
and recursively defined $c_i>0$:
$$c_i=c_{i-1}\|L(\cd)\|_qB\({\b q-1\over q-1},{\b_{i-1} q-1\over q-1}\)^{1\over q'},\qq1\les i\les k.$$
We let $k\ges1$ be the smallest integer that $\b_k\ges1$. Then the above implies
$$y(t)\les a(t)+\sum_{i=0}^{k-1}c_i\int_0^t{L(s)a(s)\over(t-s)^{1-\b_i}}ds
+c_kT^{\b_k-1}\int_0^tL(s)y(s)ds,\qq\ae~t\in[0,T].$$
Now, let
$$z(t)=\int_0^tL(s)y(s)ds,\qq t\in[0,T].$$
Then
$$\dot z(t)=L(t)y(t)\les L(t)a(t)+\sum_{i=0}^{k-1}c_iL(t)\int_0^t{L(s)a(s)\over(t-s)^{1-\b_i}}ds+c_kT^{\b_k-1}L(t)z(t).$$
Hence,
$$\ba{ll}
\ns\ds z(t)\les\int_0^te^{c_kT^{\b_k-1}\int_s^tL(\t)d\t}L(s)a(s)ds
+\sum_{i=0}^{k-1}c_{i}\int_0^te^{c_kT^{\b_k-1}\int_s^tL(\t)d\t}L(s)\int_0^s{L(\t)a(\t)
\over(s-\t)^{1-\b_i}}d\t ds\\
\ns\ds \qq\les\int_0^te^{c_kT^{\b_k-1}\int_s^tL(\t)d\t}L(s)a(s)ds
+\sum_{i=0}^{k-1}c_{i}K\int_0^tL(\t)a(\t)d\t\int_\t^t{L(s)
\over(s-\t)^{1-\b_i}}ds \\
\ns\ds\qq\les\1n\int_0^t\2n e^{c_kT^{\b_k-1}\int_s^tL(\t)d\t}L(s)a(s)ds
\1n+\2n\sum_{i=0}^{k-1}\1n c_{i}K\Big(\int_\t^t\2n L(s)^qds\Big)^{1\over q}\Big(\int_\t^t\2n{1\over(s-\t)^{(1-\b_i){q\over{q-1}}}}ds\Big)^{q-1\over q}\2n\int_0^t\2n L(s)a(s)ds \\
\ns\ds\qq\les c_k\int_0^tL(s)a(s)ds,\ea$$
for a properly redefined constant $c_k>0$. Hence,
$$y(t)\les a(t)+\sum_{i=0}^{k-1}c_i\int_0^t{L(s)a(s)\over(t-s)^{1-\b_i}}ds
+c_k\int_0^tL(s)a(s)ds,\qq\ae~t\in[0,T],$$
proving our conclusion. \endpf

\ms

Comparing with the Gronwall type inequality appearing in literature on fractional differential equations (see \cite{Henry}, for example), our inequality only involves a finite sum, instead of an infinite series.

\section{State Equation}

In this section, we discuss our state equation \eqref{1.1}, together with the cost functional \eqref{1.2}. In what follows, $U$ will be a separable metric space with the metric $\rho$, which could be a non-empty bounded or unbounded set in $\dbR^m$ with the metric induced by the usual Euclidean norm. Let $u_0\in U$ be fixed. For any $p\ges1$, we define
$$\sU^p[0,T]=\big\{u:[0,T]\to U\bigm|u(\cd)\hb{ is measurable},~\rho(u(\cd),u_0)\in L^p(0,T;\dbR)\big\}.$$

\subsection{Well-posedness in $L^p$ space}

We introduce the following assumptions for the generator $f(\cd\,,\cd\,,\cd\,,\cd)$ of our state equation.

\ms

{\bf(H1)} Let the map $f:\D\times\dbR^n\times U\to\dbR^n$ be measurable. There are nonnegative functions $L(\cd),L_0(\cd)$ with
\bel{L,L_0}L(\cd)\in L^{({1\over\b}\vee{p\over p-1})+}(0,T;\dbR),\qq
L_0(\cd)\in L^{({p\over1+\b p}\vee1)+}(0,T;\dbR),\ee
for some $p\ges1$ (with the convention that ${1\over0}=\infty$) and $\b\in(0,1)$ such that
\bel{|f-f|}\left\{\2n\ba{ll}
\ns\ds|f(t,s,y_1,u)-f(t,s,y_2,u)|\les{L(s)|y_1-y_2|\over(t-s)^{1-\b}},\qq\forall(t,s,u)\in\D\times U,~y_1,\ y_2\in\dbR^n,\\
\ns\ds|f(t,s,0,u)|\les{L(s)\rho(u,u_0)+L_0(s)\over (t-s)^{1-\b}},\qq\forall(t,s,u)\in\D\times U.\ea\right.\ee

\ms

Note that the larger the $\b\in(0,1)$, the weaker the singularity of the generator $f(\cd\,,\cd\,,\cd\,,\cd)$. Also, \eqref{|f-f|} imply
\bel{|f|*}|f(t,s,y,u)|\les{L(s)\big[|y|+\rho(u,u_0)\big]+L_0(s)\over (t-s)^{1-\b}},\qq\forall(t,s,y,u)\in\D\times\dbR^n\times U.\ee

\ms

We now present the well-posedness of the state equation \eqref{1.1} in $L^p$ spaces.

\bt{Well-posedness} \sl Let {\rm(H1)} hold with some $p\ges1$ and $\b\in(0,1)$.
Then for any $\eta(\cd)\in L^p(0,T;\dbR^n)$ and $u(\cd)\in\sU^p[0,T]$, \eqref{1.1} admits a unique solution $y(\cd)\equiv y(\cd\,;\eta(\cd),u(\cd))\in L^p(0,T;\dbR^n)$, and the following estimates hold
\bel{|y|}\|y(\cd)\|_p\les\|\eta(\cd)\|_p+K\(1+\|\rho(u(\cd),u_0)\|_p\).\ee
If $(\eta_1(\cd),u_1(\cd)),\ (\eta_2(\cd),u_2(\cd))\in L^p(0,T;\dbR^n)\times\sU^p[0,T]$ and $y_1(\cd),\ y_2(\cd)$ are the solutions of \eqref{1.1}
corresponding to $(\eta_1(\cd),u_1(\cd))$ and $(\eta_2(\cd),u_2(\cd))$, respectively, then
\bel{|y-y|}\ba{ll}
\ns\ds\|y_1(\cd)-y_2(\cd)\|_p\les K\Big\{\|\eta_1(\cd)-\eta_2(\cd)\|_p\\
\ns\ds\qq\qq\qq\qq\qq+\[\int_0^T\(\int_0^t|f(t,s,y_1(s),u_1(s))-f(t,s,y_1(s),u_2(s))|ds\)^pdt\]^{1\over p}\Big\}.\ea\ee

\et

\it Proof. \rm Fix any $\eta(\cd)\in L^p(0,T;\dbR^n)$ and $u(\cd)\in\cU^p[0,T]$. For any $z(\cd)\in L^p(0,T;\dbR^n)$, define
$$\sT[z(\cd)](t)=\eta(t)+\int_0^tf(t,s,z(s),u(s))ds,\qq t\in[0,T].$$
Denote $\th(t)=t^{\b-1}I_{(0,\infty)}(t)$, where $I_{(0,\infty)}$ is the characteristic function of $(0,\infty)$.  Then
\bel{|Tz|}\ba{ll}
\ns\ds\|\sT[z(\cd)]\|_p\les\|\eta(\cd)\|_p+\(\int_0^T\Big|\int_0^tf(t,s,z(s),u(s))ds\Big|^pdt\)^{1\over p}\\
\ns\ds\qq\qq\q\les\|\eta(\cd)\|_p+\[\int_0^T\(\int_0^t{L(s)\big[|z(s)|+\rho(u(s),u_0)\big]+L_0(s)
\over(t-s)^{1-\b}}ds\)^pdt\]^{1\over p}\\
\ns\ds\qq\qq\q\les\|\eta(\cd)\|_p+\|\th(\cd)*\big\{L(\cd)\big[\rho(u(\cd),u_0)+|z(\cd)|\big]\big\}\|_p
+\|\th(\cd)*L_0(\cd)\|_p\equiv\|\eta(\cd)\|_p+I_1+I_0.\ea\ee
Now, we split the proof into three cases.

\ms

\it Case 1. \rm $p>{1\over1-\b}$. In this case,
$${1\over\b}>{p\over p-1},\qq{p\over1+\b p}>1.$$
For any $\e\in(0,{\b\over1-\b})$, which is equivalent to $(1-\b)(1+\e)<1$, define $q$ through the following:
$${1\over q}={1\over p}+1-{1\over1+\e}<{1\over p}+1-{1\over1+{\b\over1-\b}}={1\over p}+\b<1.$$
The last inequality in the above follows from $p>{1\over1-\b}$. Thus, $1<q<p$ and
$$q\searrow{p\over1+\b p}>1,\qq\hb{as }\e\nearrow{\b\over1-\b}.$$
Since $L_0(\cd)\in L^{{p\over1+\b p}+}(0,T;\dbR)$, we
may assume that $L_0(\cd)\in L^q(0,T;\dbR)$ (for an $\e$ being close enough to
${\b\over1-\b}$). Hence, by Young's inequality,
$$I_0\equiv\|\th(\cd)*L_0(\cd)\|_p\les\|\th(\cd)\|_{1+\e}\|L_0(\cd)\|_q.$$
Also,
$${1\over q}-{1\over p}=1-{1\over1+\e}<\b.$$
Thus,
$${p-q\over pq}\nearrow\b,\qq\hb{as }\e\nearrow{\b\over1-\b},$$
which is equivalent to
$${pq\over p-q}\searrow{1\over\b},\qq\hb{as }\e\nearrow{\b\over1-\b}.$$
Hence, by $L(\cd)\in L^{{1\over\b}+}(0,T;\dbR)$, we could find $\e$ which is close enough
to ${\b\over1-\b}$ so that $L(\cd)\in L^{pq\over p-q}(0,T;\dbR)$. Then
$$\ba{ll}
\ns\ds I_1\equiv\|\th(\cd)*\big\{L(\cd)\big[|z(\cd)|+\rho(u(\cd),u_0)\big]\big\}\|_p\\
\ns\ds\q\les \|\th\|_{1+\e}\|L(\cd)\big[|z(\cd)|+\rho(u(\cd),u_0)\big]\|_q
\les\|\th(\cd)\|_{1+\e}\|L(\cd)\|_{pq\over p-q}\|\,|z(\cd)|+\rho(u(\cd),u_0)\|_p.\ea$$
Then we see that $\sT:L^p(0,T;\dbR^n)\to L^p(0,T;\dbR^n)$. Next, let $z_1(\cd),z_2(\cd)\in L^p(0,T;\dbR^n)$, we look at the following:
$$\ba{ll}
\ns\ds\|\sT[z_1(\cd)]-\sT[z_2(\cd)]\|_{L^p(0,\d;\dbR^n)}\equiv\(\int_0^\d|
\sT[z_1(\cd)](t)-\sT[z_2(\cd)](t)|^pdt\)^{1\over p}\\
\ns\ds\les\[\int_0^\d\Big|\int_0^t\(f(t,s,z_1(s),u(s))-f(t,s,z_2(s),u(s))\)ds\Big|^pdt\]^{1\over p}\\
\ns\ds\les\[\int_0^\d\(\int_0^t{L(s)|z_1(s)-z_2(s)|\over(t-s)^{1-\b}}ds\)^pdt\]^{1\over p}
\leq\|\th(\cd)*\big[L(\cd)|z_1(\cd)-z_2(\cd)|\big]\|_{L^p(0,\d;\dbR)}\\
\ns\ds\les\1n\|\th(\cd)\|_{L^{1+\e}(0,\d;\dbR)}\|L(\cd)|z_1(\cd)\1n-\1n z_2(\cd)|\,\|_{L^q(0,\d;\dbR)}\1n\les\1n\({\d^{\b\1n-\1n(1\1n-\1n\b)\e}\over\b-(1-\b)\e}\)^{1\over1+\e}\|L(\cd)\|_{L^{pq\over p-q}(0,\d;\dbR)}\|z_1(\cd)\1n-\1n z_2(\cd)\|_{L^p(0,\d;\dbR^n)}.\ea$$
Clearly, for $\d>0$ small, the map $\sT:L^p(0,\d;\dbR^n)\to L^p(0,\d;\dbR^n)$ is a contraction.
Hence, it admits a unique fixed point on $L^p(0,\d;\dbR^n)$, which is the unique solution of the state equation \eqref{1.1} on $[0,\d]$.

\ms

Next, we look \eqref{1.1} on $[\d,T]$, which can be written as
$$y(t)=\eta(t)+\int_0^\d f(t,s,y(s),u(s))ds+\int_\d^tf(t,s,y(s),u(s))ds,\qq t\in[\d,T].$$
Since (similar to \eqref{|Tz|})
$$\ba{ll}
\ns\ds\Big\|\eta(\cd)+\int_0^\d f(\cd\,,s,y(s),u(s))ds\Big\|_{L^p(\d,T;\dbR^n)}
=\[\int_\d^T\Big|\eta(t)+\int_0^\d f(t,s,y(s),u(s))ds\Big|^pdt\]^{1\over p}\\
\ns\ds\les\(\int_\d^T|\eta(t)|^pdt\)^{1\over p}+\[\int_\d^T\(\int_0^\d|f(t,s,y(s),u(s))|ds\)^pdt\]^{1\over p}\\
\ns\ds\les\|\eta(\cd)\|_{L^p(\d,T;\dbR^n)}+\[\int_\d^T\(\int_0^\d{L(s)\big[|y(s)|
+\rho(u(s),u_0)\big]+L_0(s)\over(t-s)^{1-\b}}ds\)^pdt\]^{1\over p}\\
\ns\ds\les\|\eta(\cd)\|_{L^p(\d,T;\dbR^n)}+\[\int_\d^T\(\int_0^\d{L(s)\big[|y(s)|+\rho(u(s),u_0)\big]
\over(t-s)^{1-\b}}ds\)^pdt\]^{1\over p}+\[\int_\d^T\(\int_0^\d{L_0(s)
\over(t-s)^{1-\b}}ds\)^pdt\]^{1\over p}\\
\ns\ds\les\|\eta(\cd)\|_{L^p(\d,T;\dbR^n)}+K\(\|y(\cd)\|_{L^p(0,\d;\dbR^n)}
+\|\rho(u(\cd),u_0)\|_{L^p(0,\d;\dbR)}+1\).\ea$$
Then using the same argument as above, we obtain the existence and uniqueness of the solution to the state equation on $[0,2\d]$. By induction, we could get the solvability of the state equation on $[0,T]$.

\ms

Now, let $(\eta_1(\cd),u_1(\cd)),\ (\eta_2(\cd),u_2(\cd))\in L^p(0,T;\dbR^n)\times\sU^p[0,T]$ and $y_1(\cd)$, $y_2(\cd)$ be the corresponding solutions. Then
$$\ba{ll}
\ns\ds|y_1(t)-y_2(t)|\les|\eta_1(t)-\eta_2(t)|\1n+\2n\int_0^t\2n|f(t,s,y_1(s),u_1(s))
\1n-\1n f(t,s,y_1(s),u_2(s))|ds\1n+\2n\int_0^t\2n{L(s)|y_1(s)\1n-\1n y_2(s)|\over(t-s)^{1-\b}}ds\\
\ns\ds\qq\qq\qq\equiv a(t)+\int_0^t{L(s)|y_1(s)-y_2(s)|\over(t-s)^{1-\b}}ds.\ea$$
Hence, by Lemma \ref{Gronwall},
\bel{|y-y|<}|y_1(t)-y_2(t)|\les a(t)+\sum_{i=1}^{k-1}c_i\int_0^t{L(s)a(s)\over(t-s)^{1-\b_i}}ds+c_k\int_0^tL(s)a(s)ds,\qq\ae~t\in[0,T],\ee
for some constants $c_i>0$ and $\b_i\in[\b,1)$. Consequently, similar to \eqref{|Tz|},
$$\ba{ll}
\ns\ds\|y_1(\cd)-y_2(\cd)\|_p\les K\(\int_0^Ta(t)^pdt\)^{1\over p}\\
\ns\ds\les K\Big\{\|\eta_1(\cd)-\eta_2(\cd)\|_p+\[\int_0^T\(\int_0^t|f(t,s,y_1(s),u_1(s))-f(t,s,y_1(s),u_2(s))|ds\)^pdt\]^{1\over p}\Big\},\ea$$
proving the stability estimate. We can use the similar argument to prove  this estimate to get \eqref{|y|}.

\ms

\it Case 2. \rm $1<p\les{1\over1-\b}$. In this case,
$${1\over\b}\les{p\over p-1},\qq{p\over1+\b p}\les1.$$
Also, since $1-\b\les{1\over p}<1$, for any $\e\in(0,p-1)$, the following holds:
$$1-\b\les{1\over p}<{1\over1+\e}.$$
This implies $(1-\b)(1+\e)<1$. Define $q$ through the following:
$${1\over p}<{1\over q}={1\over p}+1-{1\over1+\e}\nearrow1,\qq\hb{as }\e\nearrow p-1.$$
Then
$${1\over q}-{1\over p}=1-{1\over1+\e}\nearrow1-{1\over p}\qq\hb{as }\e\nearrow p-1.$$
Thus,
$${pq\over p-q}\searrow{p\over p-1},\qq\hb{as }\e\nearrow p-1.$$
Consequently, by choosing $\e>0$ close enough to $p-1$, we have $q>1$ close enough to 1 and ${pq\over p-1}$ close enough to ${p\over p-1}$. Hence,
$$\|\th(\cd)*L_0(\cd)\|_p\les\|\th(\cd)\|_{1+\e}\|L_0(\cd)\|_q,$$
and
$$\|\th(\cd)*\big\{L(\cd)\big[|z(\cd)|+\rho(u(\cd),u_0)\big]\big\}\|_p\les\|\th(\cd)\|_{1+\e}\|L(\cd)\|_{pq\over p-q}\|\,|z(\cd)|+\rho(u(\cd),u_0)\|_p.$$
The rest of the proof is similar to that of Case 1.

\ms

\it Case 3. \rm $p=1$. In this case, the condition reads $L(\cd)\in L^\infty(0,T;\dbR)$ and $L_0(\cd)\in L^{1+}(0,T;\dbR)$. Then
$$\|\th(\cd)*L_0(\cd)\|_1\les\|\th(\cd)\|_1\|L_0(\cd)\|_1,$$
and
$$\|\th(\cd)*\big\{L(\cd)\big[|z(\cd)|+\rho(u(\cd),u_0)\big]\big\}\|_1\les\|\th(\cd)\|_1\|L(\cd)\|_\infty
\|\,|z(\cd)|+\rho(u(\cd),u_0)\|_1.$$
The rest of the proof is similar to that of Case 1. \endpf

\ms

Let us make some comments and observations on the above theorem. First of all, the above theorem gives some sufficient conditions under which for $(\eta(\cd),u(\cd))\in L^p(0,T;\dbR^n)\times\sU^p[0,T]$, equation \eqref{1.1} admits a unique solution $y(\cd)\in L^p(0,T;\dbR^n)$. The conditions we imposed in (H1) are compatibility conditions of the integrability for the free term $\eta(\cd)$, the control $u(\cd)$, and the coefficients $L(\cd)$ and $L_0(\cd)$. From the above, we see that if $(\eta(\cd),u(\cd))\in L^p(0,T;\dbR)\times\cU^p[0,T]$ with $p>{1\over1-\b}$, then by assuming $L(\cd)\in L^{{1\over\b}+}(0,T;\dbR)$ and $L_0(\cd)\in L^{{1\over\b}-}(0,T;\dbR)$ (note that ${p\over1+\b p}<{1\over\b}$), the equation has a unique solution $y(\cd)\in L^p(0,T;\dbR^n)$. This is the case, in particular, if $\eta(\cd)\in L^\infty(0,T;\dbR^n)$ and $U$ is bounded (under the metric $\rho$). We will come back to this later. On the other hand, if $1\les p<{1\over1-\b}$, that is, say, the free term and/or the control have weaker integrability, then we need to strengthen the integrability condition for $L(\cd)$ from $L^{{1\over\b}+}$ to $L^{{p\over p-1}+}$ (in the current case, ${p\over p-1}>{1\over\b}$) to get $L^p$ solution $y(\cd)$. But, the integrability of $L_0(\cd)$ is only required to be $L^{1+}(0,T;\dbR)$. Finally, since we have used the contraction mapping theorem to establish the well-posedness of the state equation, one can see that the solution to the state equation can be obtained by a Picard iteration.

\ms

Let us present an example from which we could get some feeling about the above result.

\ms

\bex{Example 3.2} \rm Consider the following Volterra integral equation
\bel{3.7*}y(t)={1\over|t-1|^{1-\g}}+\int_0^t{\sqrt{(s-1)^{2\d-2}+y(s)^2}\over|s-1|^{1-\a}(t-s)^{1-\b}}ds,\qq\ae\ t\in[0,T],\ee
for some $\a,\b\in(0,1)$, $\g,\d\in(0,1]$, and with $T>1$. In this case, we have/can take
$$\eta(t)={1\over|t-1|^{1-\g}},\q L(s)={1\over|s-1|^{1-\a}},\q L_0(s)={1\over|s-1|^{2-\a-\d}}.$$
In order $\eta(\cd)\in L^p(0,T;\dbR)$, we need
$$p(1-\g)<1\qq\iff\qq p<{1\over1-\g}\q({1\over0}\deq\infty).$$
In order $L(\cd)\in L^{({1\over\b}\vee{p\over p-1})+}(0,T;\dbR)$, one needs
$${1-\a\over\b}<1\q\iff\q\a+\b>1,$$
and
$$(1-\a){p\over p-1}<1\q\iff\q1-\a<1-{1\over p}\q\iff\q p>{1\over\a}.$$
Finally, in order $L_0(\cd)\in L^{({p\over1+\b p}\vee1)+}(0,T;\dbR)$, one needs
$$(2-\a-\d){p\over1+\b p}<1\q\iff\q2-\a-\b-\d<{1\over p}\q\iff\q p<{1\over(2-\a-\b-\d)^+},$$
and
$$2-\a-\d<1\qq\iff\qq\a+\d>1.$$
Hence, equation \eqref{3.7*} has a unique solution $y(\cd)\in L^p(0,T;\dbR^n)$ for any $p\in[{1\over\a},{1\over(1-\g)\vee(2-\a-\b-\d)^+})$, provided
\bel{a+b>1**}\a+\b>1,\qq\a+\d>1.\ee
We point out that in general, the solution $y(\cd)$ of the equation \eqref{3.7*} is not necessarily continuous, even if the free term $\eta(\cd)$ is continuous. In fact, let $\g=1$. Then $\eta(t)\equiv1$ which is continuous. It is seen that the solution $y(\cd)$ is positive (which can be seen from a Picard iteration). Consequently,
$$\lim_{t\to1}y(t)\ges1+\lim_{t\to1}\int_0^t{ds\over|s-1|^{2-\a-\d}(t-s)^{1-\b}}
=\int_0^1{ds\over(1-s)^{3-\a-\b-\d}}=\infty,$$
provided
\bel{a+b+d<2}3-\a-\b-\d>1\qq\iff\qq\a+\b+\d<2.\ee
This will be the case if we take
$$\a={2\over3},\q\b=\d={1\over2}.$$
In this case, the solution $y(\cd)\in L^p(0,T;\dbR)$ exists with $p\in({3\over2},3)$ and it is discontinuous at $t=1$.

\ex

Note that in the above example, the solution $y(\cd)$ is discontinuous at $t=1$ only, which is the singularity of $L(\cd)$ and $L_0(\cd)$. It is natural to ask what will be the result for the general situation? Such a question has its own interest. And also since the values $y(t_i)$ of $y(\cd)$ are needed in the cost functional \eqref{1.2}, we would like to locate the discontinuity points of the solution $y(\cd)$ a priori based on the information of $L(\cd)$ and $L_0(\cd)$. This leads to the following subsection.

\subsection{Continuity of the solution}

In this subsection, we would like to explore the continuity of the solution $y(\cd)$ to the state equation \eqref{1.1}. Let us begin with some observations. Suppose $y(\cd)\in L^p(0,T;\dbR^n)$ is the unique solution to the state equation \eqref{1.1} which is rewritten here:
\bel{state*}y(t)=\eta(t)+\int_0^tf(t,s,y(s),u(s))ds,\qq t\in[0,T].\ee
Then the continuity of $y(\cd)$ is determined by that of $\eta(\cd)$ and
$$\psi(\cd)\equiv\int_0^\cd f(\cd\,,s,y(s),u(s))ds.$$
The continuity of $\eta(\cd)$ should be given a priori. Thus, we need to look at the continuity of the above-defined function $\psi(\cd)$. Hence, the preliminary results presented in Section 2 will play an interesting role here. To make it precise, we introduce the following hypothesis.

\ms

{\bf(H2)} Let $w(\cd)$ be given by \eqref{w(s)}
with  $\a_i\in(0,1)$, $0\les i\les\ell$ and $0\leq s_0<s_1<\cds<s_\ell\les T$.
Let $f:\D\times\dbR^n\times U\to\dbR^n$ be given by the following:
\bel{f}f(t,s,y,u)={f_0(t,s,y,u)\over w(s)(t-s)^{1-\b}},\qq(t,s,y,u)\in\D\times\dbR^n\times U,\ee
with $\b\in(0,1)$ and $f_0:\D\times\dbR^n\times U\to\dbR^n$ being measurable such that
\bel{f0-f0}|f_0(t,s,y,u)-f_0(t',s,y,u)|\les\o(|t-t'|),\qq\forall(t,s),\ (t',s)\in\D,~(y,u)\in\dbR^n\times U,\ee
for some modulus of continuity $\o:[0,\infty)\to[0,\infty)$, and
\bel{3.20}|f_0(t,s,y,u)|\les\bar\f(s),\qq\forall(t,s,y,u)\in\D\times\dbR^n\times U,\ee
and
\bel{3.21}|f_0(t,s,y_1,u)-f_0(t,s,y_2,u)|\les\bar L(s)(|y_1-y_2|),\qq(t,s,u)\in\D\times U,~y_1,\ y_2\in\dbR^n,\ee
for some measurable functions $\bar\f,\ \bar L:[0,T]\to[0,\infty)$.

\ms

Note that under (H2), we will have (H1) if one takes the following:
$$L_0(s)={\bar\f(s)\over w(s)},\qq L(s)={\bar L(s)\over w(s)},\qq s\in[0,T].$$
Thus, according to Theorem \ref{Well-posedness}, state equation \eqref{1.1} admits a unique solution in $L^p(0,T;\dbR^n)$, under (H2), for any $\eta(\cd)\in L^p(0,T;\dbR^n)$ with some $p\in[1,\infty)$, if
\bel{L/w}{\bar\f(\cd)\over w(\cd)}\in L^{({p\over1+\b p}\vee1)+}(0,T;\dbR),\qq{\bar L(\cd)\over w(\cd)}\in L^{({1\over\b}\vee{p\over p-1})+}(0,T;\dbR).\ee
Let us make a simple observation on the above condition. Recall the definition of $\d_0$ and $\bar s_i$ from \eqref{d0}.
It is not hard to see that \eqref{L/w} holds if $(1-\a_i)({1\over\b}\vee{p\over p-1})<1$, $0\les i\les\ell$ and for some $\e\in(0,\d_0)$, the following holds:
$$\left\{\2n\ba{ll}
\ns\ds\bar\f(\cd)\in L^{({p\over1+\b p}\vee1)+}([0,T]\setminus\bigcup_{i=0}^\ell( (s_i-\e)\vee0,(s_i+\e)\land T);\dbR),\q\\
\ns\ds\bar L(\cd)\in L^{({1\over\b}\vee{p\over p-1})+}([0,T]\setminus  \bigcup_{i=0}^\ell((s_i-\e)\vee0,(s_i+\e)\land T);\dbR),\\
\ns\ds\bar\f(\cd),\bar L(\cd)\in L^\infty((s_i-\e)\vee0,(s_i+\e)\land T;\dbR),\qq0\les i\les\ell.\ea\right.$$
Namely, due to the special structure of $w(\cd)$, it suffices to have boundedness of $\bar\f(\cd)$ and $\bar L(\cd)$ near $s_i$ $(0\les i\les\ell$) and proper integrability of these functions away from the points $s_i$. Therefore, the condition \eqref{L/w} is very mild.

\ms

We have the following result which is a direct consequence of Lemma \ref{Lemma 2.2}.

\bp{} \sl Let {\rm(H2)} hold with \eqref{L/w} for some $p\in[1,\infty)$, and $\bar\f(\cd)\in L^q(0,T;\dbR)$, $\ds q>{1\over\b}\vee{1\over\a_i}$, for all $i=0,1,\cds,\ell$. Then for any $\eta(\cd)\in L^p(0,T;\dbR^n)$ and any $u(\cd)\in\sU^p[0,T]$, state equation \eqref{1.1} admits a unique solution $y(\cd)\in L^p(0,T;\dbR^n)$ such that
\bel{wy}y(\cd)-\eta(\cd)\in L^\infty_{\bar w(\cd)}((0,T);\dbR^n)\bigcap C_{\bar w^\e(\cd)}([0,T];\dbR^n),\ee
where $\bar w(\cd)$ and $\bar w^\e(\cd)$ are given in \eqref{bar w} and \eqref{we(s)}.
\ep

\subsection{Special cases}

In this subsection, we look at some special cases.

\ms

\it 1. Linear Volterra integral equations. \rm Consider the following equation:
\bel{linear}y(t)=\eta(t)+\int_0^t{A(t,s)y(s)\over w(s)(t-s)^{1-\b}}ds,\qq t\in[0,T],\ee
where $\b\in(0,1)$, $w(\cd)$ is a weight function defined by \eqref{w(s)} and $A:\D\to\dbR^n$ satisfies %
\bel{A}|A(t,s)|\les\bar L(s),\qq\forall(t,s)\in\D,\ee
for some measurable function $\bar L(\cd)$ satisfying
\bel{barL/w}{\bar L(\cd)\over w(\cd)}\in L^{({1\over\b}\vee{p\over p-1})+}(0,T;\dbR),\ee
with some $p\in[1,\infty)$. Then, by Theorem \ref{Well-posedness}, for any $\eta(\cd)\in L^p(0,T;\dbR^n)$, equation \eqref{linear} admits a
unique solution $y(\cd)\in L^p(0,T;\dbR^n)$. Moreover, if we define operator $\cA$ by
$$\cA[y(\cd)](t)=\int_0^t{A(t,s)y(s)\over w(s)(t-s)^{1-\b}}ds,\qq t\in[0,T],$$
then, thanks to \eqref{barL/w}, by the proof of Theorem \ref{Well-posedness}, we see that $\cA:L^p(0,T;\dbR^n)\to L^p(0,T;\dbR^n)$ is a linear bounded operator. Our linear integral equation \eqref{linear} reads
$$y(\cd)=\eta(\cd)+\cA[y(\cd)].$$
Therefore, the unique solution $y(\cd)$ admits the following (abstract) representation:
$$y(\cd)=(I-\cA)^{-1}\eta(\cd)=\sum_{k=0}^\infty\cA^k\eta(\cd).$$
Now, let $(t,\t)\mapsto\F(t,\t)$ be the unique solution to the following equation:
\bel{Y}\F(t,\t)={A(t,\t)\over w(\t)(t-\t)^{1-\b}}+\int_\t^t{A(t,s)\F(s,\t)\over w(s)(t-s)^{1-\b}}ds,\qq0\les\t<t\les T.\ee
Then one has
\bel{variation of constant}y(t)=\eta(t)+\int_0^t\F(t,s)\eta(s)ds,\q\ae\ t\in[0,T].\ee
This is called the {\it variation of constant formula}.

\ms

\it 2. Fractional differential equations. \rm Let us first recall some basic notions of fractional integrals and derivatives. For $\a\in(0,1)$, let
\bel{Ia}[I^\a f(\cd)](t)={1\over\G(\a)}\int_0^t{f(s)\over(t-s)^{1-\a}}ds,\qq t\ges0,\ee
and
as long as the right hand side is well-defined, where $\G(\cd)$ is the Gamma function. We call $I^\a$ the $\a$-th order integral operator. Let
\bel{Da}[D^\a y(\cd)](t)={d\over dt}[I^{1-\a}y(\cd)](t)\equiv{1\over\G(1-\a)}{d\over dt}
\int_0^t{y(s)\over(t-s)^\a}ds,\ee
and
\bel{Da*}[D^\a_*y(\cd)](t)=[D^\a (y(\cd)-y(0))](t)=[D^\a (y(\cd)](t)-{y(0)\over\Gamma(1-\a)}t^{-\a}.\ee
In particular, when $y(\cd)\in AC([0,T];\dbR)$, the set of all absolutely continuous functions defined on $[0,T]$, one has
\bel{Da*}[D^\a_*y(\cd)](t)=[I^{1-\a}y'(\cd)](t)\equiv{1\over\G(1-\a)}\int_0^t{y'(s)\over
(t-s)^\a}ds,\ee
We call $D^\a$ and $D^\a_*$ the $\a$-th order {\it Riemann-Liouville} and {\it Caputo differential operators}, respectively. We have the following standard result (see \cite{Kilbas-Srivastava-Trujillo 2006}, Lemmas 2.5 and 2.22).

\bp{} \sl Let $\a\in(0,1)$. Then for any $y(\cd)\in L^1(0,T;\dbR)$ with $[I^{1-\a}y(\cd)](\cd)\in AC([0,T];\dbR)$.
\bel{ID}I^\a\{D^\a[y(\cd)]\}(t)=y(t)-{I^{1-\a}[y(\cd)](0)\over\G(\a)t^{1-\a}},\qq\ae~t\in(0,T];
\ee
and for $y(\cd)\in AC([0,T];\dbR)$,
\bel{ID*}I^\a\{D^\a_*[y(\cd)]\}(t)=y(t)-y(0).\ee
\ep

Now, let us consider the following fractional differential equation of Riemann-Liouville type:
\bel{Dy=f}D^\a[y(\cd)](t)=f(t,y(t),u(t)),\qq t\in[0,T].\ee
Applying the operator $I^\a$ to the above, we obtain
\bel{y=}y(t)={I^{1-\a}[y(\cd)](0)\over\G(\a)t^{1-\a}}+{1\over\G(\a)}\int_0^t{f(s,y(s),u(s))\over
(t-s)^{1-\a}}ds,\qq t\in[0,T].\ee
We refer the readers to Theorem 3.1 in \cite{Kilbas-Srivastava-Trujillo 2006} for the equivalence of (\ref{Dy=f}) and (\ref{y=}).

\ms

Likewise, if we consider the following fractional differential equation of Caputo type:
\bel{D*y=f}D^\a_*[y(\cd)](t)=f(t,y(t),u(t)),\qq t\in[0,T],\ee
applying the operator $I^\a$ to the above, we obtain
\bel{*y=}y(t)=y(0)+{1\over\G(\a)}\int_0^t{f(s,y(s),u(s))\over(t-s)^{1-\a}}ds,\qq t\in[0,T].\ee
We refer the readers to Theorem 3.24 in \cite{Kilbas-Srivastava-Trujillo 2006} for the equivalence of (\ref{D*y=f}) and (\ref{*y=}).

\ms

From the above, we see that fractional differential equations of Riemann-Liouville and Caputo types are special cases of \eqref{1.1}.

\subsection{A backward linear Volterra integral equation}

In this subsection, we consider the following linear backward Volterra integral equation:
\bel{backward}\psi(t)=\xi(t)+\int_t^T{A(s,t)^\top\psi(s)\over w(t)(s-t)^{1-\b}}ds,\qq t\in[0,T],\ee
where $A:\D\to\dbR^{n\times n}$ satisfies \eqref{A}--\eqref{barL/w}. Such an equation will play an important role in the next section. Let $1<p<{1\over1-\b}$. We claim that for any $\xi(\cd)\in L^{p\over p-1}(0,T;\dbR^n)$, the above equation admits a unique solution $\psi(\cd)\in L^{p\over p-1}(0,T;\dbR^n)$. In fact, by condition \eqref{barL/w}, we can find an $r>{1\over\b}\vee{p\over p-1}$ such that ${\bar L(\cd)\over w(\cd)}\in L^r(0,T;\dbR)$. By $r>{1\over\b}$, we can find an $\e>0$ such that
$${1\over1+\e}=1-{1\over r}>1-\b\q\Ra\q(1+\e)(1-\b)<1.$$
Then, for any $\psi(\cd)\in L^{p\over p-1}(0,T;\dbR^n)$, we have (denoting $p'={p\over p-1}$)
$$\ba{ll}
\ns\ds\Big\|\1n\int_\cd^T\2n{A(s,\cd)^\top\psi(s)\over w(\cd)(s-\cd)^{1-\b}}ds\Big\|_{p'}\2n
=\1n\Big\{\1n\int_0^T\2n\Big|\1n\int_t^T{A(s,t)^\top\2n\psi(s)\over w(t)(s-t)^{1-\b}}ds\Big|^{p'}dt\Big\}^{1\over p'}\3n\les\1n
\Big\{\1n\int_0^T\2n\({\bar L(t)\over w(t)}\int_t^T\2n{|\psi(s)|\over(s-t)^{1-\b}}ds\Big)^{p'}dt\Big\}^{1\over p'}\\
\ns\ds\les\[\int_0^T\({\bar L(t)\over w(t)}\)^rdt\]^{1\over r}\[\int_0^T\(\int_t^T{|\psi(s)|\over(s-t)^{1-\b}}ds\)^{p'r\over r-p'}dt\]^{r-p'\over p'r}\\
\ns\ds\les\Big\|{\bar L(\cd)\over w(\cd)}\Big\|_r\|\th(\cd)*|\psi(\cd)|\,\|_{p'r\over r-p'}
\les\Big\|{\bar L(\cd)\over w(\cd)}\Big\|_r\|\th(\cd)\|_{1+\e}\|\psi(\cd)\|_{p'}.\ea$$
Here, the Young's inequality for convolution is used with
$${1\over1+\e}+{1\over p'}={r-p'\over p'r}+1={1\over p'}-{1\over r}+1\q\iff\q{1\over1+\e}=1-{1\over r}.$$
By a similar argument used in the proof of Theorem \ref{Well-posedness}, we get the well-posedness of equation \eqref{backward}.

\section{Pontryagin's Maximum Principle}

In this section,  we discuss the optimal control problem for equation \eqref{1.1} with cost functional \eqref{1.2}. To begin with, let us introduce the following assumptions. The conditions assumed are more than sufficient. But for the simplicity of presentation, we prefer to use these stronger conditions.

\ms

{\bf(H3)} Let $h^j:\dbR^n\to\dbR$, $j=1,2,\cds,m$ be continuously differentiable, and  $g:[0,T]\times\dbR^n\times U\to\dbR$ be measurable with $y\mapsto g(t,y,u)$ being continuously differentiable. There exist a constant $L>0$ and a modulus of continuity $\o:[0,+\infty)\to[0,+\infty)$ such that
$$\ba{ll}
\ns\ds|g(t,y_1,u_1)-g(t,y_2,u_2)|\les L|y_1-y_2|+\o\big(\rho(u_1,u_2)\big),\qq\forall\,(t,y_1,u_1),(t,y_2,u_2)\in[0,T]\times\dbR^n\times U,\\
\ns\ds|g(t,0,u)|\les L,\qq\qq\forall(t,u)\in[0,T]\times U,\\
\ns\ds|g_y(t,y_1,u_1)-g_y(t,y_2,u_2)|\les\o(|y_1-y_2|+\rho(u_1,u_2)),\qq\forall(t,y_1,u_1),(t,y_2,u_2)\in [0,T]\times\dbR^n\times U.\ea$$
Suppose $0\les s_0<s_1<\cds<s_\ell\les T$ are given as in (H2), and $0<t_1<t_2<\cds<t_m\les T$ such that
\bel{tj}t_j\notin\{s_0,s_1,\cds,s_\ell\},\qq\forall j=1,2,\cds,m.\ee

\ms

Clearly, under (H2)--(H3), our cost functional \eqref{1.2} is well-defined. Hence, we can formulate the following optimal control problem.

\ms

{\bf Problem (P)}\ Find a $u^*(\cd)\in\sU^p[0,T]$ such that
\bel{4.1}J(u^*(\cd))=\inf_{u(\cd)\in\sU^p[0,T]}J(u(\cd)).\ee
Any $u^*(\cd)$ satisfying \eqref{4.1} is called an {\it optimal control} of Problem $(P)$, the corresponding state $y^*(\cd)$ is called an {\it optimal state} and $(y^*(\cd),u^*(\cd))$
is called an {\it optimal pair}.

\ms

In this section, we shall first give a set of necessary conditions for optimal pairs of Problem (P).
Usually, such a result is referred to as a {\it Pontryagin's maximum principle}. Then, we shall show some examples.

\subsection{ Pontryagin's maximum principle for Problem (P) }

In establishing the Pontryagin's maximum principle for the case that $U$ is not assumed to be convex, we need the following Liapunoff type theorem (see, Corollary 3.8 of Chapter 4 in \cite{Li-Yong 1995}).

\bl{LY} \sl Let $X$ be a Banach space. For any $\d>0$, let
$$\sE_\d=\big\{E\in[0,T]\bigm||E|=\d T\big\},$$
where $|E|$ stands for the Lebesgue measure of $E$. Then for any $h(\cd)\in C([0,T];L^1(0,T;X))$,
\bel{inf=0}\inf_{E\in\sE_\d}\Big\|\int_0^T\({1\over\d}{\bf1}_E(s)-1\)h(\cd,s)ds\Big\|_{C([0,T];X)}=0.\ee

\el

The following is our main result of this section, which is called Pontryagin's maximum principle for Problem $(P)$.

\bt{Theorem 4.3} \sl Let {\rm(H2)--(H3)} hold for some $p\ges1$, $\bar\f(\cd)\in L^q(0,T;\dbR)$, $\ds q>{1\over\b}\vee{1\over\a_i}$, for all $i=0,1,\cds,\ell$ and \eqref{L/w}  holds. Let $f^0:\D\times\dbR^n\times U\to\dbR$ be measurable with $y\mapsto f^0(t,s,y,u)$ being differentiable. Let $\eta(\cd)\in L^p(0,T;\dbR^n)$ and $\eta(\cd)$ be continuous at $t_j$, $j=1,2,\cds,m$. Suppose $(y^*(\cd),u^*(\cd))$ is an optimal pair of Problem {\rm(P)}. Then
there exists a solution $\psi(\cd)\in L^{p\over p-1}(0,T;\dbR^n)$ of the following {\it adjoint equation}
\bel{adjoint}\ba{ll}
\ns\ds\psi(t)=-g_y(t,y^*(t),u^*(t))^\top-\sum_{j=1}^m
{\bf1}_{[0,t_j)}(t)f_y(t_j,t,y^*(t),u^*(t))^\top h^j_y\big((y^*(t_j)\big)^\top\\
\ns\ds\qq\qq+\int_t^Tf_y(s,t,y^*(t),u^*(t))^\top\psi(s)ds,\qq t\in[0,T],\ea\ee
such that the following {\it maximum condition} holds:
\bel{max}\ba{ll}
\ns\ds\int_s^T\psi(t)^\top f(t,s,y^*(s),u^*(s))dt-g(s,y^*(s),u^*(s))-\sum_{j=1}^mh^j_y\big(y^*(t_j)\big){\bf1}_{[0,t_j]}(s)
f(t_j,s,y^*(s),u^*(s))\\
\ns\ds=\min_{u\in U}\[\int_s^T\psi(t)^\top f(t,s,y^*(s),u)dt-g(s,y^*(s),u)-\sum_{j=1}^mh^j_y\big(y^*(t_j)\big){\bf1}_{[0,t_j]}(s)
f(t_j,s,y^*(s),u)\],\\
\ns\ds\qq\qq\qq\qq\qq\qq\qq\qq\qq\qq\qq\qq\qq\qq\qq\ae~s\in[0,T].\ea\ee

\et

\it Proof. \rm We split the proof into several steps.

\ms

\it Step 1. A variational inequality. \rm Let $(y^*(\cd),u^*(\cd))$ be an optimal pair of Problem (P). Fix any $u(\cd)\in\sU^p[0,T]$. Denote
\bel{u(d)}u^\d(t)=\left\{\2n\ba{ll}
\ns\ds u^*(t),\qq t\in[0,T]\setminus E_\d,\\
\ns\ds u(t),\qq t\in E_\d,\ea\right.\ee
with $E_\d\subseteq[0,T]$ being measurable and undetermined (see {\it Step 2}). It is obvious that the control $u^\d(\cd)$ is in $\sU^p[0,T]$.
Let $y^\d(\cd)=y(\cd\,;\eta(\cd),u^\d(\cd))$ be the corresponding solution, and let
$$Y^\d(t)={y^\d(t)-y^*(t)\over\d},\qq t\in[0,T].$$
Then, $Y^\d(\cd)$ satisfies
$$\ba{ll}
\ns\ds Y^\d(t)={1\over\d}\Big\{\int_0^t\[f(t,s,y^\d(s),u^\d(s))-f(t,s,y^*(s),u^\d(s))\]ds\\
\ns\ds\qq\qq\qq+\int_0^t\[f(t,s,y^*(s),u^\d(s))-f(t,s,y^*(s),u^*(s))\]ds\Big\}\\
\ns\ds\qq\;\;=\int_0^t\[\int_0^1f_y(t,s,y^*(s)+\t\d Y^\d(s),u^\d(s))d\t\]Y^\d(s)ds\\
\ns\ds\qq\qq\qq+{1\over\d}\int_0^t{\bf1}_{E_\d}(s)\[f(t,s,y^*(s),u(s))-f(t,s,y^*(s),u^*(s))\]ds\\
\ns\ds\qq\;\;\equiv\int_0^tf^\d_y(t,s)Y^\d(s)ds+{1\over\d}\int_0^t{\bf1}_{E_\d}(s)\h f(t,s)ds,\ea$$
with
\bel{4.7}\left\{\2n\ba{ll}
\ns\ds f^\d_y(t,s)=\int_0^1f_y(t,s,y^*(s)+\t\d Y^\d(s),u^\d(s))d\t,\\ [3mm]
\ns\ds\h f(t,s)=f(t,s,y^*(s),u(s))-f(t,s,y^*(s),u^*(s)).\ea\right.\qq(t,s)\in\D.\ee
By the optimality of $(y^*(\cd),u^*(\cd))$, one has the following variational inequality:
\bel{4.20}\ba{ll}
\ns\ds0\les{J(u^\d(\cd))-J(u^*(\cd))\over\d}={1\over\d}\Big\{\int_0^T\[g(t,y^\d(t),u^\d(t))
-g(t,y^*(t),u^*(t))\]ds\\
\ns\ds\qq\qq\qq\qq\qq\qq\qq+\sum_{j=1}^m\[h^j(y^\d(t_j))-h^j(y^*(t_j))\]\Big\}\\
\ns\ds\q=\2n\int_0^T\2n\[\int_0^1g_y(t,y^*(t)+\t\d Y^\d(t),u^\d(t))d\t\]Y^\d(t)dt+{1\over\d}\int_{E_\d}\2n\big[g(t,y^*(t),u(t))-g(t,y^*(t),u^*(t))\big]dt\\
\ns\ds\qq+\sum_{j=1}^m\[\int_0^1h^j_y(y^*(t_j)+\t\d Y^\d(t_j))d\t\]Y^\d(t_j).\ea\ee

\it Step 2. Convergence of $Y^\d(\cd)$ and so on. \rm We introduce the following integral equation:
\bel{variational system}Y(t)=\int_0^t\[f_y(t,s,y^*(s),u^*(s))Y(s)+\h f(t,s)\]ds,\qq t\in[0,T],\ee
where $\h f(\cd\,,\cd)$ is given by \eqref{4.7}. Under our conditions, the above admits a unique solution $Y(\cd)$ such that
$$Y(\cd)\in L^p(0,T;\dbR^n)\bigcap\(\bigcap_{i=1}^\ell C\big((s_{i-1},s_i);\dbR^n\big)\),$$
and
$$Y(\cd)\in L^\infty_{\bar w(\cd)}(0,T;\dbR^n)\bigcap C_{\bar w^\e(\cd)}([0,T];\dbR^n).$$
We now show that for a suitable choice of $E_\d$, the following holds:
$$\lim_{\d\to0}\|Y^\d(\cd)-Y(\cd)\|_p=0,\qq\lim_{\d\to0}Y^\d(t_j)=Y(t_j),\q1\les j\les m.$$
Note that
$$\ba{ll}
\ns\ds\int_0^t{\bf1}_{E_\d}(s)\[f(t,s,y^*(s),u(s))-f(t,s,y^*(s),u^*(s))\]ds\\
\ns\ds=\int_0^t{\bf1}_{E_\d}(s){f^0(t,s,y^*(s),u(s))-f^0(t,s,y^*(s),u^*(s))\over w(s)(t-s)^{1-\b}}ds\equiv\int_0^t{\bf1}_{E_\d}(s){\h f^0(t,s)\over w(s)(t-s)^{1-\b}}ds,\ea$$
with
$$\h f^0(t,s)=f^0(t,s,y^*(s),u(s))-f^0(t,s,y^*(s),u^*(s)),\qq(t,s)\in\D.$$
By \eqref{3.20}, we have
$$|\h f^0(t,s)|\les2\bar\f(s),\qq(t,s)\in\D,$$
with $\bar\f(\cd)\in L^q(0,T;\dbR)$, where $q$ satisfies $\ds q>{1\over\b}\vee{1\over\a_i}$, for all $i=0,1,\cds,\ell$. Let
$$h(t,s)={\bf1}_{[0,t)}(s){\h f^0(t,s)\over w(s)(t-s)^{1-\b}},\qq(t,s)\in[0,T]^2.$$
Then for any $\bar p>1$ sufficiently large,
$$\ba{ll}
\ns\ds\[\int_0^T\(\int_0^T\big|h(t,s)|ds\)^{\bar p}dt\]^{1\over\bar p}
\les\[\int_0^T\(\int_0^t{2\bar\f(s)\over w(s)(t-s)^{1-\b}}ds\)^{\bar p}dt\]^{1\over\bar p}\\
\ns\ds\equiv\Big\|\th(\cd)*{2\bar\f(\cd)\over w(\cd)}\Big\|_{\bar p}\les2\|\th(\cd)\|_{1+\e}\Big\|{\bar\f(\cd)\over w(\cd)}\Big\|_q<\infty,\ea$$
where $\e>0$ is chosen so that
$${1\over1+\e}+{\,1\,\over q}={\,1\,\over\bar p}+1,\qq(1+\e)(1-\b)<1.$$
Then
$${\,1\,\over\bar p}={\,1\,\over q}+{1\over1+\e}-1=-\(\b-{\,1\,\over q}\)+{1-(1+\e)(1-\b)\over1+\e}.$$
Note that on the right hand side of the above, the first term is valued in $(-\b,0)$ since $q>{1\over\b}$; and the second term is positive with the range
$$\Big\{{1-(1+\e)(1-\b)\over1+\e}\bigm|\e\in[0,{\b\over1-\b}]\Big\}=[0,\b].$$
Thus, by suitably choosing $\e>0$, we may make $\bar p>0$ as large as we wish. Consequently, for our given $p\geq1$, by choosing $\e>0$ properly, we may have $\bar p\ges p$. Hence,
$$\[\int_0^T\(\int_0^T\big|h(t,s)|ds\)^pdt\]^{1\over p}<\infty.$$
Clearly, there exists a sequence of continuous functions $h_k(\cd\,,\cd)$ such that
$$\[\int_0^T\(\int_0^T|h(t,s)-h_k(t,s)|ds\)^pdt\]^{1\over p}<{1\over k},\qq\forall k\ges1.$$
Now, for each $h_k(\cd\,,\cd)$, applying Lemma \ref{LY}, we have that for any fixed $\d>0$, there exists some $E^k\subseteq\sE_\d$ such that
$$\sup_{t\in[0,T]}\Big|\int_0^T\({1\over\d}{\bf1}_{E^k}(s)-1\)h_k(t,s)ds\Big|<{1\over k}.$$
Then
$$\ba{ll}
\ns\ds\Big\{\int_0^T\Big|\int_0^t\({1\over\d}{\bf1}_{E^k}(s)-1\){\h f^0(t,s)\over w(s)(t-s)^{1-\b}}ds
\Big|^pdt\Big\}^{1\over p}=\Big\{\int_0^T\Big|\int_0^T\({1\over\d}{\bf1}_{E^k}(s)-1\)h(t,s)ds\Big|^pdt\Big\}^{1\over p}\\
\ns\ds\les\Big\{\int_0^T\Big|\int_0^T\({1\over\d}{\bf1}_{E^k}(s)-1\)h_k(t,s)ds\Big|^pdt\Big\}^{1\over p}\\
\ns\ds\qq+\Big\{\int_0^T\Big|\int_0^T\({1\over\d}{\bf1}_{E^k}(s)-1\)\big[h(t,s)-h_k(t,s)\big]ds\Big|^pdt\Big\}^{1\over p}\les{1\over k}\({1\over\d}+2\).\ea$$
Hence, for any fixed $\d>0$,
$$\inf_{E\in\sE_\d}\Big\{\int_0^T\Big|\int_0^t\({1\over\d}{\bf1}_E(s)-1\){\h f^0(t,s)\over w(s)(t-s)^{1-\b}}ds\Big|^pdt\Big\}^{1\over p}=0.$$
Next, for any $t_j$, let us observe
$$\int_0^T{\bf1}_{[0,t_j]}(s){|\h f^0(t_j,s)|\over w(s)(t_j-s)^{1-\b}}ds=\int_0^{t_j}{|\h f^0(t_j,s)|\over w(s)(t_j-s)^{1-\b}}ds\les\int_0^{t_j}{2\bar \f(s)\over w(s)(t_j-s)^{1-\b}}ds<\infty.$$
Hence,
$$\inf_{E\in\sE_\d}\Big|\int_0^{t_j}\({1\over\d}{\bf1}_E(s)-1\){\h f^0(t_j,s)\over w(s)(t_j-s)^{1-\b}}ds\Big|=0,\qq1\les j\les m.$$
Likewise, we also have
$$\inf_{E\in\sE_\d}\Big|\int_0^T\({1\over\d}{\bf1}_E(s)-1\)\big[g(s,y^*(s),u(s))
-g(s,y^*(s),u^*(s))\big]ds\Big|=0.$$
Now, we consider the following map
$$H(t,s)=\begin{pmatrix}\ds{\h f^0(t,s)\over w(s)(t-s)^{1-\b}}{\bf1}_{[0,t)}(s)\\
                         \ds{\h f^0(t_1,s)\over w(s)(t_1-s)^{1-\b}}{\bf1}_{[0,t_1)}(s)\\
                         \vdots\\
                         \ds{\h f^0(t_m,s)\over w(s)(t_m-s)^{1-\b}}{\bf1}_{[0,t_m)}(s)\\
                         g(s,y^*(s),u(s))-g(s,y^*(s),u^*(s))\end{pmatrix}$$
Then applying Lemma \ref{LY} to the above function in a proper product space, we obtain that for any $\d>0$, there exists an $E_\d\in\sE_\d$ such that
\bel{}\left\{\2n\ba{ll}
\ds\Big\{\int_0^T\Big|\int_0^t\({1\over\d}{\bf1}_{E_\d}(s)-1\){\h f^0(t,s)\over w(s)(t-s)^{1-\b}}ds
\Big|^pdt\Big\}^{1\over p}=o(1),\\
\ns\ds\Big|\int_0^{t_j}\({1\over\d}{\bf1}_{E_\d}(s)-1\){\h f^0(t_j,s)\over w(s)(t_j-s)^{1-\b}}ds\Big|=o(1),\qq1\les j\les m,\\
\ns\ds\Big|\int_0^T\({1\over\d}{\bf1}_{E_\d}(s)-1\)\big[g(s,y^*(s),u(s))
-g(s,y^*(s),u^*(s))\big]ds\Big|=o(1).\ea\right.\ee
By choosing such a family of $E_\d$, $\d>0$, we see that the following convergence hold:
$$\left\{\2n\ba{ll}
\ds\lim_{\d\to0}\|Y^\d(\cd)-Y(\cd)\|_p=0,\\
\ns\ds\lim_{\d\to0}Y^\d(t_j)=Y(t_j),\qq1\les j\les m,\\
\ns\ds\lim_{\d\to0}\int_0^T\({1\over\d}{\bf1}_{E_\d}(s)-1\)\big[g(s,y^*(s),u(s))
-g(s,y^*(s),u^*(s))\big]ds=0.\ea\right.$$
Hence, we end up with the following variational inequality
$$\ba{ll}
\ns\ds0\les\int_0^T\(g_y(t,y^*(t),u^*(t))Y(t)+g(t,y^*(t),u(t))-g(t,y^*(t),u^*(t))\)dt+\sum_{j=1}^m
h^j_y\big(y^*(t_j)\big)Y(t_j)\\
\ns\ds\q=\int_0^T\(g_y(s,y^*(s),u^*(s))Y(s)+g(s,y^*(s),u(s))-g(s,y^*(s),u^*(s))\)ds\\
\ns\ds\qq+\int_0^T\sum_{j=1}^mh^j_y\big(y^*(t_j)\big){\bf1}_{[0,t_j]}(s)
\(f_y(t_j,s,y^*(s),u^*(s))Y(s)+f(t_j,s,y^*(s),u(s))-f(t_j,s,y^*(s),u^*(s))\)ds\\
\ns\ds\q=\int_0^T\(g_y(s,y^*(s),u^*(s))+\sum_{j=1}^mh^j_y\big(y^*(t_j)\big){\bf1}_{[0,t_j]}
(s)f_y(t_j,s,y^*(s),u^*(s))\)Y(s)ds\\
\ns\ds\qq+\int_0^T\[g(s,y^*(s),u(s))-g(s,y^*(s),u^*(s))\\
\ns\ds\qq\qq+\sum_{j=1}^mh^j_y\big(y^*(t_j)\big){\bf1}_{[0,t_j]}(s)\(f(t_j,s,y^*(s),u(s))
-f(t_j,s,y^*(s),u^*(s))\)\]ds,\ea$$
with $Y(\cd)$ being the solution to the variational equation \eqref{variational system}.

\ms

\it Step 3. Duality. \rm Let $\psi(\cd)$ be the solution to the adjoint equation \eqref{adjoint}. Then we have
$$\ba{ll}
\ns\ds0\les\int_0^T\(g_y(s,y^*(s),u^*(s))+\sum_{j=1}^mh^j_y\big(y^*(t_j)\big){\bf1}_{[0,t_j]}
(s)f_y(t_j,s,y^*(s),u^*(s))\)Y(s)ds\\
\ns\ds\qq+\int_0^T\[g(s,y^*(s),u(s))-g(s,y^*(s),u^*(s))\\
\ns\ds\qq\qq+\sum_{j=1}^mh^j_y\big(y^*(t_j)\big){\bf1}_{[0,t_j]}(s)\(f(t_j,s,y^*(s),u(s))
-f(t_j,s,y^*(s),u^*(s))\)\]ds\\
\ns\ds=\int_0^T\(-\psi(s)+\int_s^Tf_y(t,s,y^*(s),u^*(s))^\top\psi(t)dt\)^\top Y(s)ds\\
\ns\ds\qq+\int_0^T\[g(s,y^*(s),u(s))-g(s,y^*(s),u^*(s))\\
\ns\ds\qq\qq+\sum_{j=1}^mh^j_y\big(y^*(t_j)\big){\bf1}_{[0,t_j]}(s)\(f(t_j,s,y^*(s),u(s))
-f(t_j,s,y^*(s),u^*(s))\)\]ds\\
\ns\ds=\int_0^T\psi(t)^\top\(-Y(t)+\int_0^tf_y(t,s,y^*(s),u^*(s))Y(s)ds\)dt\\
\ns\ds\qq+\int_0^T\[g(s,y^*(s),u(s))-g(s,y^*(s),u^*(s))\\
\ns\ds\qq\qq+\sum_{j=1}^mh^j_y\big(y^*(t_j)\big){\bf1}_{[0,t_j]}(s)\(f(t_j,s,y^*(s),u(s))
-f(t_j,s,y^*(s),u^*(s))\)\]ds\\
\ns\ds=\int_0^T\[-\psi(t)^\top\int_0^t\(f(t,s,y^*(s),u(s))-f(t,s,y^*(s),u^*(s))\)ds\]dt\\
\ns\ds\qq\qq+\int_0^T\[g(s,y^*(s),u(s))-g(s,y^*(s),u^*(s))\\
\ns\ds\qq\qq+\sum_{j=1}^mh^j_y\big(y^*(t_j)\big){\bf1}_{[0,t_j]}(s)\(f(t_j,s,y^*(s),u(s))
-f(t_j,s,y^*(s),u^*(s))\)\]ds\\
\ns\ds=\int_0^T\[\int_s^T-\psi(t)^\top\(f(t,s,y^*(s),u(s))-f(t,s,y^*(s),u^*(s))\)dt\\
\ns\ds\qq\qq+g(s,y^*(s),u(s))-g(s,y^*(s),u^*(s))\\
\ns\ds\qq\qq+\sum_{j=1}^mh^j_y\big(y^*(t_j)\big){\bf1}_{[0,t_j]}(s)\(f(t_j,s,y^*(s),u(s))
-f(t_j,s,y^*(s),u^*(s))\)\]ds.\ea$$
Hence, using the Lebesgue point theorem for integrable functions, we reach the following:
$$\ba{ll}
\ns\ds\int_s^T\psi(t)^\top f(t,s,y^*(s),u^*(s))dt-g(s,y^*(s),u^*(s))-\sum_{j=1}^mh^j_y\big(y^*(t_j)\big){\bf1}_{[0,t_j]}(s)
f(t_j,s,y^*(s),u^*(s))\\
\ns\ds\ges\[\int_s^T\psi(t)^\top f(t,s,y^*(s),u)dt-g(s,y^*(s),u)-\sum_{j=1}^mh^j_y\big(y^*(t_j)\big){\bf1}_{[0,t_j]}(s)
f(t_j,s,y^*(s),u)\],\ \ae~s\in[0,T].\ea$$
This gives the maximum condition \eqref{max}. \endpf

\subsection{ Special cases in the sense of Riemann-Liouville and Caputo senses}

In recent years, optimal control problems for fractional differential equations have attracted the attention of some researchers. However, most of the works on maximum principles for fractional differential equations were
established by convex perturbation technique. See, for instance, Agrawal \cite{Agrawal 2004}, Agrawal--Defterli--Baleanu \cite{Agrawal-Defterli-Baleanu 2010}, Frederico--Torres \cite{Frederico-Torres 2008} and Kamocki \cite{Kamocki 2014b}
in the sense of Riemann-Liouville case, and  Agrawal  \cite{Agrawal 2008}, Bourdin \cite{Bourdin 2012} and   Hasan--Tangpong--Agrawal \cite{Hasan-Tangpong-Agrawal} in the sense of Caputo case.

\ms

Let us take a look a recent work \cite{Kamocki 2014b}, in which Kamocki considered the fractional differential equation of Riemann-Liouville type (\ref{Dy=f}) with $\a\in(0,1)$ and some convex assumptions. The control $u(\cd)$ takes value in a compact set $U$ in $\dbR^m$ and $f$ satisfies
$$\ba{ll}
\ns\ds|f(t,y_1,u)-f(t,y_2,u)|\les N|y_1-y_2|,\qq\forall\,y_1,y_2\in \dbR^n,\,\ t\in[0,T], \ u\in U,\\
\ns\ds|f(t,0,u)|\les r(t)+\g|u|,\qq\qq\forall(t,u)\in[0,T]\times U,\ea$$
where $N>0$ and $\g\ges 0$ are two constants and $r(\cd)\in L^p(0,T;\dbR)$.
The corresponding solution belongs to $L^p(0,T;\dbR^n)$ for some $p\ges1$.
When  $p>1$ and $I^{1-\a}[y(\cd)](0)=0$, a Pontryagin's maximum principle for Problem (P) was proved. For the case $I^{1-\a}[y(\cd)](0)\ne0$,
maximum principle was obtained only for $1<p<\frac{1}{1-\a}$.

\ms

It is easy to check that all the above-mentioned results for fractional differential equations are the special cases of what we presented in the previous subsection.

\section{Concluding Remarks}

This paper presented some analysis of singular Volterra integral equations, and established a Pontryagin type maximum principle for an optimal control of such kind of equations. Here are some remarks in order.

\ms

$\bullet$ As we have indicated, the fractional differential equations of Riemann-Liouville or Caputo types of order no more than one are fully covered by our results. For fractional differential equations of higher order, similar results can be obtained by properly modifying our approach.

\ms

$\bullet$ It is easy to see that all the results that we presented will remain true for non-singular Volterra integral equations.

\ms

$\bullet$ We have allowed to have very general singularity in the free term and the generator. Therefore, our results can apply to a much wider class of problems than those covered by fractional differential equations and non-singular Volterra integral equations.

\vfil\eject


\begin{thebibliography}{1}{\small}



\rm

\bibitem{Agrawal 2004} O.~P.~Agrawal, \it A general formulation and solution scheme for fractional optimal control problems, \rm Nonlinear Dynamics, \rm 38 (2004), 323-337.

\bibitem{Agrawal 2008}O.~P.~Agrawal, \it A formulation and numerical scheme for fractional optimal control problems, \rm Journal of Vibration and Control, \rm 14(9-10) (2008), 1291--1299.


\bibitem{Agrawal-Defterli-Baleanu 2010} O.~P.~Agrawal, O.~Defterli, and D.~Baleanu, \it Fractional optimal control problems with several state and control variables, \sl J. Vibration \& Control,
   \rm 16 (2010), 1967--1976.


\bibitem{Angell 1976} T.~S.~Angell, \it On the optimal control of systems governed by nonlinear integral equations, \sl J. Optim. Theory Appl., \rm 19 (1976), 63--79.

\bibitem{Arafa-Rida-Khalil 2012} A.~A.~M.~Arafa, S.~Z.~Rida, and M.~Khalil, \it Solutions of fractional order model of childhook diseases with constant vaccination strategy, \sl Math. Sci. Lett., \rm 1 (2012), 17--23.


\bibitem{Belbas 2007} S.~A.~Belbas, \it A new method for optimal control of Volterra integral equations, \sl Applied Math. Comp., \rm 189 (2007), \rm 1902--1915.

\bibitem{Belbas 2008} S.~A.~Belbas, \it A reduction method for optimal control of Volterra integral equations, \sl Appl. Math. Comp., \rm 197 (2008), 880--890.

\bibitem{Benson 1998} D.~A.~Benson, \it The fractional advection-dispersion equation: development and application, \sl Ph.D. thesis, \rm University of Nevada at Reno, 1998.

\bibitem{Bonnans-de la Vega-Dupuis 2013} J.~F.~Bonnans, C.~de la Vega, and X.~Dupuis, \it First- and second-order optimality conditions for optimal control problems of state constrained integral equations, \sl J. Optim. Theory Appl., \rm 159 (2013), 1--40.

\bibitem{Bogachev 2007} V.~I.~Bogachev, \sl Measure Theory, I, \rm Springer-Verlag, New York, 2007.

\bibitem{Bourdin 2012} L.~Bourdin, \it A class of fractional optimal control problems and fractional Pontryagin's systems. Existence of a fractional Noether's theorem, \rm arXive:1203.1422v1, 2012.

\bibitem{Burnap-Kazemi 1999} C.~Burnap and M.~A.~Kazemi, \it Optimal control of a system governed by nonlinear Volterra integral equations with delay, \sl IMA J. Math. Control Inform., \rm 16 (1999), 73--89.

\bibitem{Caputo 1967} M.~Caputo, \sl Linear models of dissipation whose $Q$ is almost frequency independent-II, \sl Geophys. J. Roy. Astron. Soc., \rm 13 (1967), 529--539; reprinted in \sl Fract. Calc. Appl. Anal., \rm 11 (2008), 4--14.

\bibitem{Caputo--Mainardi 1971} M.~Caputo and F.~Mainardi, \it A new dissipation model based on memory mechanism, \sl Pure Appl. Geophys., \rm 91 (1971), 134--147; reprinted in \sl Fract. Calc. Appl. Anal., \rm 10 (2007), 310--323.


\bibitem{Carlson 1987} D.~A.~Carlson, \it An elementary proof of the maximum principle for optimal control problems governed by a Volterra integral equation, \sl J. Optim. Theory Appl., \rm 54 (1987), 43--61.

\bibitem{Chern 1993} J.~T.~Chern, \it Finite element modeling of viscoelastic materials on the theory of fractional calculus, \sl Ph.D. thesis, \rm Pennsylvania State University, 1993.

\bibitem{Dalir-Bashour 2010} M.~Dalir and M.~Bashour, \it Applications of fractional calculus, \sl Appl. Math. Sci., \rm 4 (2010), 1021--1032.

\bibitem{Das-Gupta 2011} S.~Das and P.~K.~Gupta, \it A mathematical model on fractional Lotka-Volterra equations, \sl J. Theoretical Biology, \rm 277 (2011), 1--6.

\bibitem{de la Vega 2006} C.~de la Vega, \it Necessary conditions for optimal terminal time control problems governed by a Volterra integral equation, \sl J. Optim. Theory Appl., \rm 130 (2006), 79--93.

\bibitem{Demirci-Unal-Ozalp 2011} E.~Demirci, A.~Unal, and N.~\"Ozalp, \it A fractional order SEIR model with density dependent death rate, \sl Hacettepe J. Math. Stat., \rm 40 (2011), 287--295.

\bibitem{Diethelm 2013} K.~Diethelm, \it A fractional calculus based model for the simulation of an outbreak of dengue fever, \sl Nonlinear Dynamics, \rm 71 (2013), 613--619.

\bibitem{Diethelm-Freed 1999} K.~Diethelm and A.~D.~Freed, \it On the solution of nonlinear fractional differential equations used in the modeling of viscoplasticity, {\rm in:}  F. Keil, W. Mackens, H. Vo$\beta$, J. Werther (eds.), \sl Scientific Computing in Chemical Engineering II: Computational Fluid Dynamics, Reaction Engineering, and Molecular Properties, \rm 217--224, Springer, Heidelberg 1999.


\bibitem{Frederico-Torres 2008} G.~S.~F.~Frederico and D.~F.~M.~Torres, \it Fractional conservation laws in optimal control theory, Nonlinear Dynamics, \rm 53 (2008), 215--222.





\bibitem{Gomez-Aguilar-Razo-Hernandez-Granados-Lieberman 2014} J.~F.~G\'omez-Aguilar, R.~Razo-Hernandez, and D.~Granados-Lieberman, \it A physical interpretation of fractional calculus in observables terms: analysis of the fractional time constant and the transitory response, \sl Revista Mexicana de Fisica, \rm 60 (2014), 32--38.

\bibitem{Hasan-Tangpong-Agrawal} M.~ M.~Hasan, X.~W.~Tangpong, and O.~P.~Agrawal,
   \it Fractional optimal control of distributed systems in spherical and cylindrical coordinates, \sl Journal of Vibration and Control, \rm 18(10)
   (2011), 1506--1525.


\bibitem{Henry} D.~Henry, \sl Geometric Theory of Semilinear Parabolic Equations, \rm Springer-Verlag, New York, 1981.


\bibitem{Kamien-Muller 1976} M.~I.~Kamien and E.~Muller, \it Optimal control with integral state equations, \sl Rev. Econ. Stud., \rm 43 (1976), 469--473.

\bibitem{Kamocki 2014} R.~Kamocki, \it On the existence of optimal solutions to fractional optimal
control problems, \sl Appl. Math. \& Computation, \rm 235 (2014), 94--104.

\bibitem{Kamocki 2014b} R.~Kamocki, \it Pontryagin maximum principle for fractional ordinary optimal control problems, \rm Math. Meth. Appl. Sci., 37 (2014), 1668--1686.

\bibitem{Kisela-2008} T.~Kisela, \it Fractional differential equations and their applications, \rm Diploma Thesis, BRNO Univ. of Technology, 2008.

\bibitem{Kilbas-Srivastava-Trujillo 2006} A.~A.~Kilbas, H.~M.~Srivastava, and J.~J.~Trujillo, \sl Theory and Applications of Fractional Differential Eqautions, \rm North-Holland, Amsterdam, 2006.

\bibitem{Li-Yong 1995} X.~Li and J.~Yong, \sl Optimal Control Theory for Infinite Dimensional Systems, \rm Birkh\"auser, Boston, 1995.

\bibitem{Medhin 1986} N.~G.~Medhin, \it Optimal process governed by integral equations, \sl J. Math. Anal. Appl., \rm 120 (1986), 1--12.

\bibitem{Metzler-Schick-Kilian-Nonnenmacher 1995} R.~Metzler, W.~Schick, H.~G.~Kilian, and T.~F.~Nonnenmacher, \it Relaxation in filled polymers: a fractional calculus approach, \sl J. Chem. Phys., \rm 103 (1995), 7180--7186.

\bibitem{Moshrefi-Torbati--Hammond 1998} M.~Moshrefi-Torbati and J.~K.~Hommond, \it Physical and geometrical interpretation of fractional operators, \sl J. Franklin Inst., \rm 335B (1998), 1077--1086.

\bibitem{Okyere-Oduro-Ampojsah-Dontwi-Frempong 2016} E.~Okyere, F.~T.~Oduro, S.~K.~Amponsah, I.~K.~Dontwi, and N.~K.~Frempng, \it Fractional order SIR model with constant population, \sl British J. Math. \& Computer Sci., \rm 14(2) (2016), 1--12.

\bibitem{Oldham-Spanier 1974} K.~B.~Oldham, J.~Spanier, \sl The Fractional Calculus, \rm Academic Press, New York, 1974.



\bibitem{Rahimy 2010} M.~Rahimy, \it Applications of fractional differential equations, \sl Appl. Math. Sci., \rm 4 (2010), no.50, 2453--2461.

\bibitem{Samko-Kilbas-Marichev-1987} S.~G.~Samko, A.~A.~Kilbas, and O.~I.~Marichev, \sl Fractional Intgrals and Applications, Theory and Applications, \rm Gordon and Breach Science Publishers, 1987.

\bibitem{Scalas-Gorenflo-Mainardi 2000} E.~Scalas, R.~Gorenflo, and F.~Mainardi, \it Fractional calculus and continuous-time finance, \sl Physica A, \rm 284 (2000), 376--384.

\bibitem{Scalas-Gorenflo-Mainardi 2004} E.~Scalas, R.~Gorenflo, and F.~Mainardi, \it Uncoupled continuous-time random walks: analytic solution and limiting behaviour of the master equation, \sl Phys. Rev. E, \rm 69 (2004), 011107.

\bibitem{Tarasov 2013} V.~E.~Tarasov, \it Review of some promising fractional physical models, \sl Int. J. Modern Physics B, \rm 27 (2013), 1330005.

\bibitem{Tenreiro Machado et al 2010} J.~A.~Tenreiro Machado, M.~F.~Silva, R.~S.~Barbosa, I.~S.~Jesus, C.~ M.~Reis, M.~G.~Marcos, and A.~F.~Galhano, \it Some applications of fractional calculus in engineering, \sl Math. Prob. Engineering, \rm 2010, Article ID 639801.

\bibitem{Torvik-Bagley 1984} P.~J.~Torvik and R.~L.~Bagley, \it On the appearance of the fractional derivative in the behavior of real materials, \sl J. Appl. Mech., \rm 51 (1984), 294--298.

\bibitem{Vinokurov 1969} V.~R.~Vinokurov, \it Optimal control of processes describted by integra equations, parts I, II, and III, \sl SIAM J. Control, \rm 7 (1969), 324--355; comments by L.~W.~Neustadt and J.~Warga, \sl SIAM J. Control, \rm 8 (1970), 572.

\bibitem{Zhang-Li-Chen 1989} X.~Zhang, X.~Li, and Z.~Chen, \sl Differential Equation Theory of Optimal Control Systems (in Chinese), \rm Higher Education Press, Beijing, 1989.


\end{thebibliography}
\end{document}